\documentclass[11pt]{amsart}
\usepackage{amssymb}
\usepackage[curve,matrix,arrow]{xy}
\DeclareMathOperator{\Id}{Id}
\DeclareMathOperator{\colim}{colim}

\DeclareMathOperator{\Hom}{Hom}
\DeclareMathOperator{\map}{map}
\DeclareMathOperator{\Tor}{Tor}
\DeclareMathOperator{\sgn}{sgn}
\DeclareMathOperator{\sh}{sh}
\DeclareMathOperator{\sk}{sk}

\newcommand{\Ab}{{\mathcal A}b}
\newcommand{\mN}{{\mathbb N}}
\newcommand{\mQ}{{\mathbb Q}}
\newcommand{\mZ}{{\mathbb Z}}
\newcommand{\Pc}{{\mathcal P}}
\newcommand{\iso}{\cong}
\newcommand{\sm}{\wedge}
\newcommand{\tensor}{\otimes}

\renewcommand{\to}{\longrightarrow}
\newcommand{\un}{\underline}
\newcommand{\xra}{\xrightarrow}
\newcommand{\spec}{Sp^{\Sigma}}

\numberwithin{equation}{section}

\newtheorem{lemma}[equation]{Lemma}
\newtheorem{prop}[equation]{Proposition}
\newtheorem{cor}[equation]{Corollary}
\theoremstyle{definition}
\newtheorem{defn}[equation]{Definition}
\newtheorem{construction}[equation]{Construction}
\newtheorem{rk}[equation]{Remark}
\newtheorem{eg}[equation]{Example}

\begin{document}

\title[Homotopy groups of symmetric spectra]
{On the homotopy groups of symmetric spectra}

\date{\today; 2000 AMS Math.\ Subj.\ Class.: 55P42, 55U35}
\author{Stefan Schwede}
\address{Mathematisches Institut, Universit\"at Bonn, Germany}
\email{schwede@math.uni-bonn.de}
\maketitle

Symmetric spectra are an easy-to-define and convenient model 
for the stable homotopy category with a nice smash product.
Symmetric ring spectra first showed up under the name `FSP on spheres`
in the context of algebraic $K$-theory and topological Hochschild homology.
Around~1993, Jeff Smith made the crucial observation that the 
`FSP on spheres` are the monoids in a category of `symmetric spectra`
with respect to an associative and commutative smash prodct, 
and he suspected compatible model category structures 
so that one obtains as homotopy categories 
`the` stable homotopy category (for symmetric spectra),
the homotopy category of $A_\infty$ ring spectra (for symmetric ring spectra),
respectively the homotopy category of $E_\infty$ ring spectra 
(for commutative symmetric ring spectra).
The details of various model structures were worked out by Hovey, Shipley and
Smith in~\cite{HSS}.

Maybe the only tricky point with symmetric spectra is that the
stable equivalences can {\bf not} be defined by looking at 
stable homotopy groups (defined as the classical sequential
colimit of the unstable homotopy groups of the terms in a
symmetric spectrum).
Formally inverting the {\em $\pi_*$-isomorphisms},
i.e., those morphisms which induce isomorphisms of stable
homotopy groups, leaves too many homotopy types.
Instead, Hovey, Shipley and Smith introduce a strictly larger
class of {\em stable equivalences}, defined as the morphisms which induce
isomorphisms on all cohomology theories.
The difference between $\pi_*$-isomorphisms and stable equivalences
has previously confused at least the present author.
The precise relationship between the naively defined homotopy groups
and the `true` homotopy groups (morphisms from sphere in the
stable homotopy category) have largely been mysterious
(although Shipley's detection functor~\cite[Sec.~3]{shipley-THH}
sheds considerable light on this).

In this paper we advertise and systematically exploit 
extra algebraic structure 
on the (classical) homotopy groups of a symmetric spectrum which,
in the authors opinion, clarifies several otherwise adhoc observations 
and illuminates various mysterious points in the theory of symmetric spectra.
This extra structure is an action of the monoid $M$
of injective self-maps of the set of natural numbers.
The $M$-modules that come up, however, have a special property
which we call {\em tameness}, see Definition~\ref{defn-tame}.
Tameness has strong algebraic consequences and severely restricts 
the kinds of $M$-modules which can arise as
homotopy groups of symmetric spectra.

Here is a first example of the use of the $M$-action.
An important class of symmetric spectra is formed by the
{\em semistable} symmetric spectra. Within this class, 
stable equivalences coincide with $\pi_*$-isomorphisms, so it is very
useful to recognize a given symmetric spectrum as semistable. 
In Proposition~\ref{prop-semistable}, we characterize
the semistable symmetric spectra as those for which the
$M$-action on homotopy groups is trivial. 
Other examples for the use of the $M$-action are:
\begin{itemize}
\item the group $\pi_0(F_nS^n)$ is not just an infinitely generated
free abelian group, but a very prominent tame $M$-module $\Pc_n$;
these $M$-modules are pairwise non-iso\-morphic 
for different~$n$ (Example~\ref{eg-free}).
\item certain monomorphisms on homotopy groups which appear
in the discussion of semistable symmetric spectra are given by
the action of the special element $d$ of $M$ defined
by $d(i)=i+1$ (Example~\ref{eg-d acts}).
\item a tame $M$-module which is finitely generated as an abelian group
necessarily has trivial $M$-action 
(Lemma~\ref{lemma-algebraic consequences}~(iv)); 
this implies that symmetric spectra
whose homotopy groups are dimensionwise finitely generated are always
semistable. 
This puts~\cite[Prop.~5.6.4~(1)]{HSS} into perspective,
which proves that dimensionwise {\em finite} homotopy groups implies
semistability.
\item the homotopy groups of the symmetric spectra $\Omega X$ and
$S^1\sm X$ (Example~\ref{eg-loop and suspension}), 
$\sh X$ (Example~\ref{eg-shift}), 
$F_1S^0\sm X$ (Example~\ref{eg-drift}), 
$RX$ and $R^\infty X$ (Example~\ref{eg-R infty})
are functors of the homotopy groups of $X$, and the $M$-action
on $\pi_*X$ determines the $M$-action on the homotopy groups of these
constructions.
\item the $E^2$-term of Shipley's Bousfield-Kan spectral sequence 
for calculating the homotopy groups of $DX$ can be identified
with the Tor groups of $\pi_*X$ over the monoid ring of $M$, 
and the spectral sequence can be analyzed in many examples
(Section~\ref{sec-spectral sequence}).
 \end{itemize}

The $M$-action is also intertwined in an interesting and non-trivial way
with the smash product of symmetric spectra.
This explains the curious phenomenon that the homotopy groups 
of a symmetric ring spectrum do in general
{\em not} form a graded ring in any natural way (unless the underlying
symmetric spectrum in semistable); the problem is that the
termwise smash product pairings are not compatible with the stabilization
maps in the colimit defining stable homotopy groups.
The natural algebraic structure on the homotopy of a symmetric ring spectrum
is an action by the injection operad of the set of natural numbers. 
This operad is a discrete analog, with similar properties, 
of the operad of linear isometries which features
prominently in other foundational parts of stable homotopy theory.
We do not discuss these topics here but hope to return to them elsewhere.

{\bf Acknowledgments.} 
We warn the reader that this paper does not contain any
difficult mathematics. We combine elementary observations centered
around the tame $M$-action on the homotopy groups of symmetric spectra
and revisit many familiar examples and phenomena from this point of view.

Moreover, several of the things which we discuss 
are contained, one way or another, explicitly or
implicitly, in the papers~\cite{HSS} and~\cite{shipley-THH}.
The $M$-action makes a brief appearance in 
Proposition~2.2.9 of~\cite{shipley-THH}, which Shipley attributes to Smith;
it appears though that the $M$-action has not been used systematically, 
nor has tameness been exploited before.

\newpage

\section{Homotopy groups as $M$-modules}

A {\em symmetric spectrum} consists of the following data:
\begin{itemize}
\item a sequence of pointed spaces
$X_n$ for $n\geq 0$;
\item a base-point preserving continuous left action of
the symmetric group $\Sigma_n$ on $X_n$, for each $n\geq 0$; 
\item pointed continuous maps $\sigma_n:X_n\sm S^1\to X_{n+1}$ for $n\geq 0$.
\end{itemize}
This data is subject to the following condition: the composite
$$\xymatrix@C=15mm{
X_n\sm S^m \ar^-{\sigma_n\,\sm\,\Id}[r] &
X_{n+1}\sm S^{m-1} \ar^-{\sigma_{n+1}\sm \Id}[r] & 
\quad \cdots \quad \ar^-{\sigma_{n+m-1}}[r] & \ X_{n+m} 
}$$
is $\Sigma_n\times\Sigma_m$-equivariant for all $n,m\geq 0$. Here $S^m$ is the
$m$-fold smash product of $S^1$, on which the symmetric group 
$\Sigma_m$ acts by permuting the factors, and the action on the
target $X_{n+m}$ is by restriction from $\Sigma_{n+m}$ to the
subgroup $\Sigma_n\times\Sigma_m$;
here and in the rest of the paper we
view $\Sigma_n\times\Sigma_m$ as a subgroup of $\Sigma_{n+m}$
via the monomorphism which sends
$(\tau,\delta)\in\Sigma_n\times\Sigma_m$ to 
$\tau\times\delta\in\Sigma_{n+m}$ defined by
$$ (\tau\times\delta)(i) \ = \ \begin{cases}
\qquad \tau(i) & \text{\quad for $1\leq i\leq n$, and}\\
\delta(i-n)+n & \text{\quad for $n+1\leq i\leq n+m$.}
\end{cases}$$

A {\em morphism} $f:X\to Y$ of symmetric spectra consists of
$\Sigma_n$-equivariant continuous pointed maps $f_n:X_n\to Y_n$ for $n\geq 0$,
which are compatible with the structure maps in the sense that
$f_{n+1}\circ\sigma_n=\sigma_n\circ (f_n\sm \Id)$
for all $n\geq 0$.
The category of symmetric spectra is denoted by $\spec$.

The {\em $k$-th stable homotopy group} of a symmetric spectrum $X$,
for $k$ any integer, is the colimit
\[ \pi_k X \ = \ \colim_n \, \pi_{k+n} X_n \ , \]
taken over the maps 
$\iota_*:\pi_{k+n} X_n \to \pi_{k+n+1} X_{n+1}$ defined as the composite
\begin{equation}\label{colimit system pi_k}
\pi_{k+n} \, X_n \ \xrightarrow{\ -\sm S^1\ } \
\pi_{k+n+1} \, \left( X_n\sm S^1 \right) \ \xrightarrow{\ (\sigma_n)_*\ } \
\pi_{k+n+1} \, X_{n+1} \ . \end{equation}

A {\em $\pi_*$-isomorphism} is a morphism of symmetric spectra
which induces an isomorphism on all stable homotopy groups. 
Every $\pi_*$-isomorphism is a stable equivalence in the sense 
of~\cite[Def.~3.1.3]{HSS}, but not conversely.
Inverting only the $\pi_*$-isomorphisms of symmetric spectra
would leave too many stable homotopy types, and the resulting category
could not be equivalent to the usual stable homotopy category.

The definition of homotopy groups does not take the
symmetric group actions into account; using these actions we will now
see that the homotopy groups of a symmetric spectrum have more structure.

\begin{construction}\label{con-M action}
We define an action of the monoid $M$ of injective self-maps of 
the set $\omega=\{1,2,3,\dots\}$ of positive natural numbers, 
on the homotopy groups of a symmetric spectrum $X$. 
We break the construction up into two steps and
pass through the intermediate category of {\em $I$-functors}.
The category $I$ has an object $\mathbf n=\{1,\dots,n\}$
for every non-negative integer $n$, including $\mathbf{0}=\emptyset$.
Morphisms in $I$ all injective maps.
An {\em $I$-functor} is a covariant functor from the category $I$
to the category of abelian groups.

{\em Step 1: from symmetric spectra to $I$-functors.}
For every integer $k$ we assign an $I$-functor $\un{\pi}_kX$ 
to the symmetric spectrum $X$. On objects, this $I$-functor is given by
$$ (\un{\pi}_kX)(\mathbf n) \ = \ \pi_{k+n}X_n $$
if $k+n\geq 2$ and $(\un{\pi}_kX)(\mathbf n)=0$ for $k+n<2$.
If $\alpha:\mathbf n\to \mathbf m$ is an injective map and $k+n\geq 2$,
then $\alpha_*: (\un{\pi}_kX)(\mathbf n)\to (\un{\pi}_kX)(\mathbf m)$
is given as follows. We choose a permutation $\gamma\in\Sigma_m$
such that $\gamma(i)=\alpha(i)$ for all $i=1,\dots,n$ and set
$$ \alpha_*(x)\ =\  \sgn(\gamma) \cdot \gamma_*(\iota_*^{m-n}(x)) $$
where $\iota_*:\pi_{k+n}X_n\to\pi_{k+n+1}X_{n+1}$ is the 
stabilization map~\eqref{colimit system pi_k}.

We have to justify that this definition is independent 
of the choice of permutation~$\gamma$. Suppose 
$\gamma'\in\Sigma_m$ is another permutation which agrees with $\alpha$
on $\mathbf n$. Then $\gamma^{-1}\gamma'$ is a permutation
of $\mathbf m$ which fixed the numbers $1,\dots,n$, so it is of the form
$\gamma^{-1}\gamma'=1\times\tau$ for some $\tau\in\Sigma_{m-n}$, 
where $1$ is the unit of $\Sigma_n$.
It suffices to show that for such permutations 
the induced action  on $\pi_{k+m}X_m$ 
via the action on $X_m$ satisfies the relation
\begin{equation} \label{action relation}
(1\times\tau)_*(\iota_*^{m-n}(x)) \ = \ \sgn(\tau)\cdot (\iota_*^{m-n}(x)) 
\end{equation}
for all $x\in\pi_{k+n}X_n$.
To justify this we let $f:S^{k+n}\to X_n$ represent $x$.
Since the iterated structure map
$\sigma^{m-n}:X_n\sm S^{m-n}\to X_m$ 
is $\Sigma_{n}\times\Sigma_{m-n}$-equivariant,
we have a commutative diagram
$$\xymatrix@C=15mm{ S^{k+m} \ar[r]^-{f\sm \Id} \ar[d]_{\Id\sm\tau} &
X_n\sm S^{m-n} \ar[r]^-{\sigma^{m-n}} \ar[d]^{\Id\sm\tau} &
X_m \ar[d]^{1\times\tau}\\
 S^{k+m} \ar[r]_-{f\sm\Id} &
X_n\sm S^{m-n} \ar[r]_-{\sigma^{m-n}} & X_m }$$
The composite through the upper 
right corner represents $(1\times\tau)_*(\iota_*^{m-n}(x))$.
Since the effect on homotopy groups of precomposing with a 
coordinate permutation of the sphere is
multiplication by the sign of the permutation,
the composite through the lower left 
corner represents $\sgn(\tau)\cdot(\iota_*^{m-n}(x))$.
This proves formula~\eqref{action relation} and completes the definition
of $\alpha_*: (\un{\pi}_kX)(\mathbf n)\to (\un{\pi}_kX)(\mathbf m)$.

The inclusion $\mathbf n\to \mathbf{n+1}$ induces the map $\iota_*$ 
over which the colimit $\pi_kX$ is formed, so if we denote the inclusion
by $\iota$, then two meanings of $\iota_*$ are consistent. 
We let $\mN$ denote the subcategory of $I$ which contains all objects 
but only the inclusions as morphisms, and then we have
$$ \pi_kX \ = \ \colim_{\mN}\ \un{\pi}_k X \ . $$

{\em Step 2: from $I$-functors to tame $M$-modules.}
The next observation is that for any $I$-functor $F$ the colimit
of $F$, formed over the subcategory $\mN$ of inclusions, 
has a natural left action by the monoid $M$ of injective self-maps
of the set $\omega$ of natural numbers. Applied to the $I$-functor
$\un{\pi}_kX$ coming from a symmetric spectrum $X$, this yields
the $M$-action on the stable homotopy group $\pi_kX$.

We let $I_\omega$ denote the category with objects the sets 
$\mathbf n$ for $n\geq 0$ and the set $\omega$ and with all injective maps
as morphisms. So $I_\omega$ contains $I$ as a full subcategory and 
contains one more object~$\omega$ whose endomorphism monoid is~$M$.
We will now extend an $I$-functor $F$ to a functor from the category $I_\omega$
in such a way that the value of
the extension at the object $\omega$ is the colimit of $F$, 
formed over the subcategory $\mN$ of inclusions.
It will thus be convenient, and suggestive, to denote the colimit of $F$, 
formed over the subcategory $\mN$ of inclusions, by $F(\omega)$
and not introduce new notation for the extended functor.
The $M$-action on the colimit of $F$ is then the action of the
endomorphisms of $\omega$ in $I_\omega$ on $F(\omega)$.

So we set $F(\omega)=\colim_\mN F$ and first define
$\beta_*:F(\mathbf n)\to F(\omega)$ for every injection 
$\beta:\mathbf n\to\omega$
as follows. We set $m=\max\{\beta(\mathbf n)\}$, 
denote by $\beta|_{\mathbf n}:\mathbf n\to \mathbf m$ 
the restriction of $\beta$ and take $\beta_*(x)$ 
to be the class in the colimit 
represented by the image of $x$  under 
$$ (\beta|_{\mathbf n})_*\ :\  F(\mathbf n)\ \to \ F(\mathbf m) \ . $$
It is straightforward to check that this is a functorial extension of
$F$, i.e., for every morphism $\alpha:\mathbf k\to\mathbf n$ in $I$ we have
$(\beta\alpha)_*(x)=\beta_*(\alpha_*(x))$.

Now we let $f:\omega\to\omega$ be an injective self-map of $\omega$,
and we want to define $f_*:F(\omega)\to F(\omega)$.
If $[x]\in F(\omega)$ is an element in the colimit represented 
by a class $x\in F(\mathbf n)$, then we set
$f_*[x]=[(f|_{\mathbf n})_*(x)]$ where $f|_{\mathbf n}:\mathbf n\to \omega$
is the restriction of $f$ and $f_*:F(\mathbf n)\to F(\omega)$
was defined in the previous paragraph. Again it is straightforward to
check that this definition does not depend in the representative $x$
of the class $[x]$ in the colimit and that the extension is functorial,
i.e., we have
$(f\alpha)_*(x)=f_*(\alpha_*(x))$ for injections 
$\alpha:\mathbf n\to \omega$ as well as
$(fg)_*[x]=f_*(g_*[x])$ when $g$ is another injective self-map of $\omega$.
As an example, if we also write $\iota:\mathbf n\to\omega$ for the inclusion,
then we have $\iota_*(x)=[x]$ for $x\in F(\mathbf n)$.

The definition just given is in fact the universal way to extend
an $I$-functor $F$ to a functor on the category $I_\omega$,
i.e., we have just constructed a left Kan extension of $F:I\to \Ab$
along the inclusion $I\to I_\omega$. However, we do not need this fact,
so we omit the proof.
\end{construction}

A trivial but important observation straight from the definition
is that the action of the monoid $M$ on the colimit of
any $I$-functor $F$ has a special property: 
every element in the colimit $F(\omega)$ is represented by a class
$x\in F(\mathbf n)$ for some $n\geq 0$; 
then for every element $f\in M$ which fixes 
the numbers $1,\dots,n$, we have $f_*[x]=[x]$.
We introduce a special name for such $M$-modules:

\begin{defn}\label{defn-tame} An {\em $M$-module} is a left
module over the monoid ring $\mZ[M]$ of the monoid $M$ of
injective self-maps of the set $\omega=\{1,2,3,\dots\}$.
We call an $M$-module $W$ {\em tame} 
if for every element $x\in W$ there exists a number 
$n\geq 0$ with the following property:
for every element $f\in M$ which fixes the set $\mathbf n$ elementwise
we have $fx=x$.
\end{defn}

Thus the homotopy groups of a symmetric spectrum are always tame $M$-modules.
An example of an $M$-module which is not tame is the free module
of rank~1.
Tameness has many algebraic consequences which we discuss in the next section. 
We show in Remark~\ref{rk-all operations} that we have now 
found all natural operations on the homotopy groups of a symmetric spectrum;
more precisely, we show that the ring of natural operations on $\pi_0X$
is a completion of the monoid ring of~$M$, so that an arbitrary operation
is a sum, possibly infinite, of operations by elements from $M$.

The action of the monoid $M$ in used by Shipley 
in~\cite[Prop.~2.2.9]{shipley-THH},
who credits this observation to Jeff Smith. To our knowledge,
the tameness of this action, which is elementary but crucial 
for many things which we do in this paper, has not been exploited before.

\begin{eg}\label{eg-d acts} 
To illustrate the action of the monoid $M$ on the homotopy groups
of a symmetric spectrum $X$ we make it explicit for the injection
$d:\omega\to\omega$ given by $d(i)=i+1$, which will also play an important role
later. For every $n\geq 1$, the map $d$ and the cycle $(1,2,\dots,n,n+1)$
of $\Sigma_{n+1}$ agree on $\mathbf n$, so $d$ acts on $\pi_k X$ as the colimit
of the system
$$\xymatrix@C=12mm{ \pi_k X_0 \ar[r]^-{\iota_*}\ar[d]_{\iota_*} & 
\pi_{k+1} X_1 \ar[r]^-{\iota_*}\ar[d]^{-(1,2)_*\circ {\iota_*}} &
\pi_{k+2} X_2 \ar[r]^-{\iota_*}\ar[d]^{(1,2,3)_*\circ {\iota_*}} & \cdots \ar[r]^-{\iota_*}& 
\pi_{k+n} X_n \ar[r]^-{\iota_*} \ar[d]^{(-1)^n(1,2,\dots,n,n+1)_*\circ \iota_*}&\\
 \pi_{k+1} X_1 \ar[r]_-{\iota_*} & \pi_{k+2} X_2 \ar[r]_-{\iota_*} &
\pi_{k+3} X_3 \ar[r]_-{\iota_*} &\cdots \ar[r]_-{\iota_*} &
\pi_{k+n+1} X_{n+1} \ar[r]_-{\iota_*}&} $$
(at least for $k\geq 0$; for negative values of $k$ only
a later portion of the system makes sense).
\end{eg}

\begin{rk}\label{rk-stable vs unstable} The stable homotopy group
$\pi_kX$ of a symmetric spectrum $X$ can also be calculated from the
system of {\em stable} as opposed to {\em unstable} homotopy groups
of the individual spaces $X_n$. For us, the $m$th stable homotopy group 
$\pi_m^{\text{s}}K$
of a pointed space~$K$ is the colimit of the sequence of abelian groups
\begin{equation}\label{eq-stable homotopy colimit} 
\pi_mK \ \xra{S^1\sm \ } \ \pi_{1+m}(S^1\sm K) \ 
\xra{S^1\sm \ } \ \pi_{2+m}(S^2\sm K) \ \xra{S^1\sm \ }\ \cdots 
\end{equation}
where we stabilize from the left. Smashing with the identity of $S^1$ 
from the right provides a stabilization map
(even an isomorphism) $\pi_m^{\text{s}}K\to\pi_{m+1}^{\text{s}}(K\sm S^1)$.
For a symmetric spectrum we can then define an 
$I$-functor $\un{\pi}_k^{\text{s}} X$ by setting 
$(\un{\pi}_k^{\text{s}} X)(\mathbf n)=\pi_{k+n}^{\text{s}} X_n$
on objects (with no restriction on $k+n$) and defining the action of 
a morphism $\mathbf n\to\mathbf m$ in the same way as for the
$I$-functor $\un{\pi}_k X$ of unstable homotopy groups.

The map from the initial term to the colimit 
of the sequence~\eqref{eq-stable homotopy colimit} 
provides a natural transformation 
$\pi_mK\to\pi_m^{\text{s}}K$ which is compatible with stabilization,
so it defines a morphism of $I$-functors 
$\un{\pi}_kX\to \un{\pi}_k^{\text{s}} X$ for every symmetric spectrum~$X$. 
The induced map on colimits
$$ \colim_\mN \ \un{\pi}_k X \ \xra{\ \iso\ } \ 
\colim_\mN \ \un{\pi}^{\text{s}}_k X$$
is bijective (compare~\cite[Lemma 2.2.3]{shipley-THH}), 
thus an isomorphism of $M$-modules.
\end{rk}

\begin{rk} We have chosen to let the spheres act from the right on
the spaces in a symmetric spectrum, and not from the left
as in~\cite{HSS}. The reason for this is that we find the $M$-action 
on the homotopy groups more transparent this way.
More formally, we consider {\em  right} $S$-modules 
in the category of symmetric
sequences where Hovey, Shipley and Smith consider {\em left} modules.
Since the sphere spectrum is a commutative monoid in symmetric sequences,
its categories of left and right modules are isomorphic.
However, the isomorphism is slightly more subtle than one may first think,
and making consistent identifications of homotopy groups can be 
a little tricky.  
The reader is invited to translate our definition of the $M$-action
on the homotopy groups of a right $S$-module via these isomorphism
to left $S$-modules.
After this translation, the action of $d\in M$ as in Example~\ref{eg-d acts}
becomes the map that occurs in the proof of Lemma~5.6.3 of~\cite{HSS}.
\end{rk}

\section{Algebraic properties of tame $M$-modules}

In this section we discuss some algebraic  properties of tame $M$-modules. 
It turns out that tameness is a rather restrictive condition.
An example is that a tame $M$-module whose underlying abelian group 
is finitely generated must necessarily have a trivial $M$-action.

We introduce some useful notation and terminology.
For an injective map $f:\omega\to \omega$ we write $|f|$ for the smallest 
number $i\geq 0$ such that $f(i+1)\ne i+1$. So in particular,
$f$ restricts to the identity on $\{1,\dots, |f|\}$. We write $|\Id|=\infty$.
An element $x$ of an $M$-module $W$ has {\em filtration $n$} 
if for every $f\in M$ with $|f|\geq n$
we have $fx=x$. We denote by $W^{(n)}$ the subgroup of $W$ of elements of
filtration $n$; for example, $W^{(0)}$ is the set of
elements fixed by all $f\in M$. We say that $x$ has 
{\em filtration exactly $n$} if it lies
in $W^{(n)}$  but not in $W^{(n-1)}$. By definition, an $M$-module
$W$ is tame if and only if every element has a finite filtration,
i.e., if the groups $W^{(n)}$ exhaust $W$. 

The following lemmas collect some elementary observations,
first for arbitrary $M$-modules and then for tame $M$-modules.

\begin{lemma}\label{lemma-not nec tame properties} Let $W$ be any $M$-module.
\begin{itemize}
\item[(i)] If two elements $f$ and $g$ of $M$ coincide on 
$\mathbf n=\{1,\dots,n\}$, then $fx=gx$ for all $x\in W$ of filtration~$n$.
\item[(ii)] For $n\geq 0$ and $f\in M$ set $m=\max\{f(\mathbf n)\}$.
Then $f\cdot W^{(n)}\subseteq W^{(m)}$.
\item[(iii)] We denote by $d\in M$ the map given by $d(i)=i+1$.
If $x\in W$ has  filtration exactly~$n$ with $n\geq 1$, then
$dx$ has filtration exactly $n+1$.
\item[(iv)] Let $V\subseteq W$ be an $M$-submodule such that
  the action of $M$ on $V$ and $W/V$ is trivial. Then the action
  of $M$ on $W$ is also trivial.
\end{itemize}
\end{lemma}
\begin{proof}
(i) We can choose a bijection $\gamma\in M$
which agrees with $f$ and $g$ on $\mathbf n$, and then $\gamma^{-1}f$ and
$\gamma^{-1}g$ fix $\mathbf n$ elementwise. So for $x$ of
filtration~$n$ we have $(\gamma^{-1}f)x=x=(\gamma^{-1}g)x$. 
Multiplying by $\gamma$ gives $fx=gx$.

(ii) If $g\in M$ satisfies $|g|\geq m$, then $gf$ and $f$ 
agree on $\mathbf n$. So for all $x\in W^{(n)}$ we have $gfx=fx$ by~(i),
which proves that $fx\in W^{(m)}$. 

(iii) We have $d\cdot W^{(n)}\subseteq W^{(n+1)}$ by part~(ii).
To prove that $d$ increases the exact filtration we consider $x\in W^{(n)}$
with $n\geq 1$ and show that $dx\in W^{(n)}$ implies $x\in W^{(n-1)}$.

For $f\in M$ with $|f|=n-1$ we define $g\in M$ by $g(1)=1$ and
$g(i)=f(i-1)+1$ for $i\geq 2$. Then we have $gd=df$ and $|g|=n$.
We let $h$ be the cycle $h=(f(n)+1, f(n),\dots,2,1)$ so that we have
$|hd|=f(n)=\max\{f(\mathbf n)\}$.
Then $fx\in W^{(f(n))}$ by part~(ii) and so
$$ fx\ =\ (hd)(fx)\ =\ h(g(dx))\ =\ (hd)x \ = \ x \ . $$
Altogether this proves that $x\in W^{(n-1)}$.

(iv) Since the $M$-action is trivial on $V$ and $W/V$,
every $f\in M$ determines an additive map $\delta_f:W/V\to V$ such that
$x-fx=\delta_f(x+V)$ for all $x\in W$. These maps satisfy 
$\delta_{fg}(x)=\delta_f(x)+\delta_g(x)$ and so $\delta$ is a homomorphism from
the monoid $M$ to the abelian group of additive maps from $W/V$ to $V$.
We will see in Lemma~\ref{BM contractible} below that the classifying space
$BM$ of the monoid $M$ is contractible, so $H^1(BM,A)=\Hom(M,A)$
is trivial for every abelian group $A$. Thus $\delta_f=0$ for all
$f\in M$, i.e., $M$ acts trivially on~$W$.
\end{proof}

\begin{cor} The assignment $\mathbf n\mapsto W^{(n)}$ extends to an
$I$-functor $W^{(\bullet)}$ in such a way that $W\mapsto W^{(\bullet)}$
is right adjoint to the functor which assigns to an $I$-functor $F$
the $M$-module $F(\omega)$. The counit of the adjunction
$(W^{(\bullet)})(\omega)\to W$ is injective with image 
the subgroup of elements of finite filtration, 
which is also the largest tame submodule of $W$. 
The assignment $W\mapsto (W^{(\bullet)})(\omega)=\bigcup_n W^{(n)}$ 
is right adjoint to the inclusion of tame $M$-modules into all $M$-modules.
\end{cor}
\begin{proof}
To define the $I$-functor
$W^{(\bullet)}$ on morphisms $\alpha:\mathbf n\to\mathbf m$ in the category $I$
we choose any extension $\tilde\alpha:\omega\to\omega$ of
$\alpha$ and define $\alpha_*:W^{(n)}\to W^{(m)}$ as the restriction
of $\tilde\alpha\cdot:W\to W$. This really has image in $W^{(m)}$ by
part~(ii) of Lemma~\ref{lemma-not nec tame properties} 
and is independent of the extension by~(i) of that lemma.
The rest is immediate.
\end{proof}

\begin{lemma}\label{lemma-algebraic consequences} Let $W$ be a tame $M$-module.
\begin{itemize}
\item[(i)] Every element of $M$ acts injectively on $W$.
\item[(ii)] If the filtration of elements of $W$ is bounded, 
then $W$ is a trivial $M$-module.
\item[(iii)] If the map $d$ given by $d(i)=i+1$ acts surjectively on $W$, 
then $W$ is a trivial $M$-module.
\item[(iv)] If $W$ is finitely generated as an abelian group, then
$W$ is a trivial $M$-module.
\end{itemize}
\end{lemma}
\begin{proof} (i) Consider $f\in M$ and $x\in W^{(n)}$ with $fx=0$.
Since $f$ is injective, we can choose $h\in M$ with $|hf|\geq n$.
Then $x=(hf)x=h(fx)=0$, so $f$ acts injectively.

(ii) Lemma~\ref{lemma-not nec tame properties}~(iii) implies that if $W=W^{(n)}$
for some $n\geq 0$, then $n=0$, so the $M$-action is trivial.

(iii) Suppose $M$ does not act trivially, so that $W^{(0)}\ne W$.
Let $n$ be the smallest positive integer such that $W^{(0)}\ne W^{(n)}$.
Then by part (iii) of Lemma~\ref{lemma-not nec tame properties}, 
any $x\in W^{(n)}-W^{(0)}$ is not in the image of $d$,
so $d$ does not act surjectively.

(iv) The union of the nested sequence of subgroups 
$W^{(0)}\subseteq W^{(1)}\subseteq W^{(2)}\subseteq\cdots$ is $W$. 
Since finitely generated abelian groups are noetherian,
we have $W^{(n)}=W$ for all large enough $n$. By part (ii), the monoid
$M$ must act trivially.
\end{proof}

Parts~(i), (iii) and~(iv) of Lemma~\ref{lemma-algebraic consequences}
can fail for non-tame $M$-modules: we can let $f\in M$ act on 
the abelian group $\mZ$ as the identity if the image of $f:\omega\to\omega$
has finite complement, and we let $f$ acts as~0 if its image has infinite
complement.

\begin{eg} \label{ex-P_n}
We introduce certain tame $M$-modules $\Pc_n$ for $n\geq 0$
which play important roles throughout this paper.
The module $\Pc_n$ is the free abelian group with basis
the set of ordered $n$-tuples of pairwise distinct elements of $\omega$
(or equivalently the set of injective maps from 
$\mathbf n$ to $\omega$).
The monoid $M$ acts from the left on this basis by componentwise evaluation, 
i.e., $f(x_1,\dots,x_n)=(f(x_1),\dots,f(x_n))$, and it acts
on $\Pc_n$ by additive extension.
For $n=0$, there is only one basis element, the empty tuple, and so
$\Pc_0$ is isomorphic to $\mZ$ with trivial $M$-action.
For $n\geq 1$, the basis is countably infinite and the $M$-action is
non-trivial. The module $\Pc_n$ is tame: the filtration of a
basis element $(x_1,\dots,x_n)$ is the maximum of the components. 
So the filtration subgroup $\Pc_n^{(m)}$ is generated by the $n$-tuples 
all of whose components are less than or equal to $m$.
An equivalent way of saying this is that 
$\Pc_n^{(m)}=\mZ[I(\mathbf n,\mathbf m)]$, the free abelian group
generated by all injections from $\mathbf n$ to $\mathbf m$;
in particular, $\Pc_n^{(m)}$ is trivial for $m<n$.

The module $\Pc_n$ represents the functor of
taking elements of filtration~$n$: for every $M$-module $W$, the map
$$ \Hom_{M\text{-mod}}(\Pc_n, W) \ \to \ W^{(n)} \ , \quad
\varphi \mapsto \varphi(1,\dots,n)$$
is bijective.
\end{eg}

\section{Examples}\label{sec-example}

We discuss several classes of symmetric spectra with a view
towards the $M$-action on the stable homotopy groups.

\begin{eg}[Stabilizing homotopy groups]\label{eg-stabilizing} 
Let $X$ be a symmetric spectrum whose homotopy groups
stabilize, i.e., for each $k\in\mZ$
there exists an $n\geq 0$ such that from the group $\pi_{k+n}X_n$
on, all maps in the system~\eqref{colimit system pi_k}
defining $\pi_kX$ are isomorphisms. 
We claim that then the action of $M$ is trivial.

To prove the claim we consider more generally any $I$-functor $F$ 
for which there exists an $n\geq 0$ such that the natural map
$F(\mathbf n)\to F(\omega)$ to the colimit is surjective.
This is certainly the case if from $\mathbf n$ 
on all maps $\iota_*$ in the colimit system are isomorphisms.
In this situation every element of $F(\mathbf n)$ has filtration~$n$. 
But tame $M$-modules with bounded filtration 
necessarily have trivial $M$-action
by Lemma~\ref{lemma-algebraic consequences}~(ii).

Examples of symmetric spectra with stabilizing homotopy groups
include all suspension spectra, $\Omega$-spectra, or $\Omega$-spectra
from some point $X_n$ on.
The symmetric spectrum obtained from a $\Gamma$-space $A$~\cite{segal}
by evaluation on spheres is another example since 
the structure map $A(S^n)\sm S^1\to A(S^{n+1})$
is $(2n+1)$-connected~\cite[Prop.~5.21]{lydakis}.
\end{eg}

\begin{eg}[Orthogonal spectra] \label{M acts trivially on orthogonal spectra}
The monoid $M$ acts trivially on the homotopy group of every symmetric
spectrum which is underlying an orthogonal spectrum~\cite[Ex.~4.4]{MMSS}.
Indeed, the inclusion $\Sigma_n\to O(n)$ as permutation matrices
sends all even permutations to the path component of the unit in $O(n)$.
So if the $\Sigma_n$-action on a pointed space $X_n$ extends to an
$O(n)$-action, then all even permutations act as the identity on
the homotopy groups of $X_n$.

So we consider more generally any $I$-functor $F$ which takes all
even permutations to identity maps.
Given $f\in M$ and an element $[x]\in F(\omega)$ in the colimit 
represented by $x\in F(\mathbf n)$
we can find $m\geq \max\{f(\mathbf n)\}$ 
and an {\em even} permutation $\gamma\in\Sigma_m$ such that
$\gamma$ agrees with $f$ on ${\mathbf n}$.
Since $\gamma$ is even, we then have 
$f_*[x]=[(f|_{\mathbf n})_*(x)]=[(\gamma|_{\mathbf n})_*(x)]
=[\gamma_*(\iota^{m-n}_*(x))]=[\iota^{m-n}_*(x)]=[x]$.
\end{eg}

\begin{eg}[Eilenberg-Mac Lane spectra]
\label{ex-EM ring spectra with M}
Every tame $M$-module $W$ can be realized as the homotopy group
of a symmetric spectrum. For this purpose we modify the construction
of the symmetric Eilenberg-Mac Lane spectrum of an abelian group,
cf.~\cite[Ex.~1.2.5]{HSS}.
We define
\[ (HW)_n \ = \ | W^{(n)}\tensor \mZ[S^n] | \ , \]
where $W^{(n)}$ is the  filtration~$n$ subgroup of $W$,
$\mZ[S^n]$ refers to the simplicial abelian group freely generated
by the simplicial set $S^n=S^1\sm\dots\sm S^1$, 
divided by the subgroup generated by the
basepoint, and the bars denote geometric realization.
The symmetric group $\Sigma_n$ takes $W^{(n)}$ to itself and we let it
act diagonally on $(HW)_n$, i.e., on $S^n$ by permuting the smash factors.
The homotopy groups of the symmetric spectrum $HW$ are concentrated
in dimension zero where we have $\pi_0 HW\iso W$ as $M$-modules. 
If $M$ acts trivially on $W$, then this is just
the ordinary Eilenberg-Mac Lane spectrum.
(Instead of the system $\mathbf n\mapsto W^{(n)}$ we could use any 
$I$-functor in the definition above; this shows that every
$I$-functor arises as the $I$-functor $\un{\pi}_0$ of a symmetric spectrum).
\end{eg}

\begin{eg}[Free symmetric spectra]\label{eg-free}
Hovey, Shipley and Smith observe in~\cite[Example~3.1.10]{HSS} 
that the zeroth stable homotopy group
of the free symmetric spectrum $F_nS^n$ is free abelian of
countably infinite rank for $n\geq 1$. We refine this calculation
to an isomorphism of $M$-modules $\pi_0(F_nS^n)\iso\Pc_n$;
here $\Pc_n$ is the $M$-module which represents taking filtration~$n$
elements, see Example~\ref{ex-P_n}. 
So while the groups  $\pi_0(F_nS^n)$ are all additively isomorphic
for different positive $n$, the $M$-action distinguishes them.
In particular, there cannot be a chain of
$\pi_*$-isomorphisms between $F_nS^n$ and $F_mS^m$ for $n\ne m$.

As in~\cite[Def.~2.2.5]{HSS} we denote by~$F_nK$ the
{\em free symmetric spectrum} generated by a pointed space $K$ in level $n$.
The free functor $F_n$ is left adjoint to evaluating a symmetric spectrum
at level~$n$. 

We claim that there is a natural isomorphism of $M$-modules
\begin{equation}\label{homotopy of free spectrum} 
\pi_k(F_nK) \ \iso \ \Pc_n\tensor\pi_{k+n}^{\text{s}}K \ . 
\end{equation}
Here  $\pi_{k+n}^{\text{s}}K$ is the $(k+n)$th stable homotopy group of $K$;
the monoid $M$ acts only on~$\Pc_n$. 
We postpone the proof to the next paragraph, 
where this isomorphism will be a special case of a more general statement.
\end{eg}

\begin{eg}[Semifree symmetric spectra]\label{ex-semifree}
Let $L$ be a pointed space with a basepoint preserving left action by the
symmetric group $\Sigma_n$, for some $n\geq 0$. 
We let $H_nL$ denote the {\em semifree symmetric spectrum} 
generated by $L$, which we define as follows. First we let $L[n]$ 
denote the symmetric sequence which has the $\Sigma_n$-space $L$ 
in level $n$ and a point everywhere else, and then we set
$H_nL=S\tensor L[n]$.
The functor  $H_n$ is left adjoint to evaluating a symmetric spectrum 
at level~$n$, but now viewed as a functor with values 
in pointed $\Sigma_n$-spaces.

We construct a natural isomorphism of $M$-modules
\begin{equation}\label{homotopy of semifree spectrum} 
\pi_k (H_n L) \ \iso \  \Pc_n\tensor_{\Sigma_n}(\pi_{k+n}^{\text{s}}L)(\sgn)  \ . 
\end{equation}
Here we use the right $\Sigma_n$-action on the tame $M$-module $\Pc_n$ 
given on the basis by permuting the components of an $n$-tuple, i.e., 
$(x_1,\dots,x_n)\gamma=(x_{\gamma(1)},\dots,x_{\gamma(n)})$.
On the right of the tensor symbol, the group $\Sigma_n$ acts by what is
induced on stable homotopy groups by the action on $L$, twisted by sign. 
For a pointed space $K$ (without any group action)
we have $H_n(\Sigma_n^+\sm K)\iso F_nK$,
which makes the isomorphism~\eqref{homotopy of free spectrum} 
a special case of~\eqref{homotopy of semifree spectrum}.

The construction of the isomorphism~\eqref{homotopy of semifree spectrum}
starts from the more explicit description of the semifree spectrum as
\[ (H_nL)_{n+m} \ = \ 
\Sigma_{n+m}^+\sm_{\Sigma_n\times \Sigma_{m}} (L\sm S^{m})  \]
($H_nL$ consists only of a point in levels less than $n$).
The structure map 
$$ \sigma_{n+m}\ : \
\left(\Sigma_{n+m}^+\sm_{\Sigma_n\times \Sigma_m} (L\sm S^m)\right) \sm S^1
\ \to \ \Sigma_{n+m+1}^+\sm_{\Sigma_n\times \Sigma_{m+1}} (L\sm S^{m+1}) $$
arises from the inclusion $\Sigma_{n+m}\to\Sigma_{n+m+1}$ (as the subgroup
fixing the element $n+m+1$) and the identification 
$S^m\sm S^1\iso S^{m+1}$.

We calculate the $I$-functor $\un{\pi}_k^{\text{s}}(H_nL)$
consisting of the {\em stable} homotopy groups of the spaces
$(H_nL)_{k+n}$ and exploit that for any symmetric spectrum $X$ the 
$M$-module $\pi_kX$ can also be calculated as the colimit of
$\un{\pi}_k^{\text{s}}X$ instead of the $I$-functor 
$\un{\pi}_kX$ of {\em unstable} homotopy groups,
see Remark~\ref{rk-stable vs unstable}.
We start with the isomorphism of $\Sigma_{n+m}$-modules
\begin{align*} \pi_{k+n+m}^{\text{s}}(H_nL)_{n+m} \ &= \ 
\pi_{k+n+m}^{\text{s}} \left(\Sigma_{n+m}^+\sm_{\Sigma_n\times \Sigma_{m}}
(L\sm S^{m})\right) \\
&\iso \ 
\mZ[\Sigma_{n+m}]\tensor_{\mZ[\Sigma_n\times\Sigma_{m}]}
(\pi_{k+n}^{\text{s}}L\tensor\sgn_m) \ . \end{align*}
The coordinate permutations of $S^m$ act by sign on stable homotopy groups,
whence the sign representation $\sgn_m$ of $\Sigma_m$.
For any $\Sigma_n$-module $B$ the map
\begin{align*} \mZ[\Sigma_{n+m}]
\tensor_{\mZ[\Sigma_n\times \Sigma_m]} (B\tensor\sgn_m) \ &\to \ 
\Pc_n^{(n+m)}\tensor_{\Sigma_n} B(\sgn) \\
\gamma\tensor (b\tensor 1) \hspace*{2cm} &\longmapsto \  
\sgn(\gamma)\cdot(\gamma(1),\dots,\gamma(n))\tensor b
\end{align*}
is an isomorphism of $\Sigma_{n+m}$-modules, where
$\Pc_n^{(n+m)}$ is the filtration $n+m$ subgroup of the $M$-module $\Pc_n$. 
Taking $B=\pi_{k+n}^{\text{s}}L$ and combining all the above yields
isomorphisms of $\Sigma_{n+m}$-modules
$$\pi_{k+n+m}^{\text{s}}(H_nL)_{n+m} \ \iso \
\Pc_n^{(n+m)}\tensor_{\Sigma_n}(\pi_{k+n}^{\text{s}}L)(\sgn) $$
which as $n+m$ varies constitute an isomorphism of $I$-functors
$$ \un{\pi}_{k}^{\text{s}} (H_nL) \ \iso \
\Pc_n^{(-)}\tensor_{\Sigma_n}(\pi_{k+n}^{\text{s}}L)(\sgn) \ . $$
Taking colimits gives the isomorphism of 
$M$-modules~\eqref{homotopy of semifree spectrum}.
\end{eg}

\begin{eg}[Infinite products]\label{ex-products}
Finite products of symmetric spectra are $\pi_*$-iso\-morphic to finite wedges, 
so stable homotopy groups commute with finite products. 
But homotopy groups do not in general commute with infinite products.
This should not be surprising because stable homotopy groups involves
a sequential colimit, and these generally do not preserve infinite products.

There are even two different ways in which commutation with products can fail.
First we note that an infinite product of a family $\{W_i\}_{i\in I}$
of tame $M$-modules is only tame if almost all the modules $W_i$
have trivial $M$-action.
Indeed, if there are infinitely many $W_i$ with non-trivial $M$-action,
then by Lemma~\ref{lemma-algebraic consequences}~(ii)
the product $\prod_{i\in I}W_i$ contains tuples of elements
whose filtrations are not bounded. We define the {\em tame product} of the
family $\{W_i\}_{i\in I}$ by
$$ \prod_{i\in I}^{\text{tame}}\, W_i \ = \ 
\bigcup_{n\geq 0} \left( \prod_{i\in I} W_i^{(n)} \right) \ , $$
which is the largest tame submodule of the product and thus
the categorical product in the category of tame $M$-modules.

Now we consider a family $\{X_i\}_{i\in I}$ of symmetric spectra.
Since the monoid $M$ acts tamely on the homotopy groups 
of any symmetric spectrum, 
the natural map from the homotopy groups of the product
spectrum to the product of the homotopy groups always lands in the
tame product. But in general, this natural map
\begin{equation}\label{product map}
\pi_k\left(\prod_{i\in I}X_i\right) \ \to \ 
\prod^{\text{tame}}_{i\in I} \, \pi_k X_i \end{equation}
need not be an isomorphism.
As an example we consider the symmetric spectra $(F_1S^1)^{\leq i}$ 
obtained by truncating 
the free symmetric spectrum $F_1S^1$ above level~$i$, i.e., 
$$ ((F_1S^1)^{\leq i})_n \ = \ \begin{cases}
(F_1S^1)_n & \text{\ for $n\leq i$,}\\
\quad *    & \text{\ for $n\geq i+1$}
\end{cases}$$
with structure maps as a quotient spectrum of $F_1S^1$.
Then $(F_1S^1)^{\leq i}$ has trivial homotopy groups for all $i$. 
The 0th homotopy group of the product $\prod_{i\geq 1}(F_1S^1)^{\leq i}$ 
is the colimit of the sequence of maps
$$  \prod_{i\geq n} \Pc_1^{(n)} \ \to \ \prod_{i\geq n+1} \Pc_1^{(n+1)} 
$$
which first projects away from the factor indexed by $i=n$ and then
takes a product of inclusions $\Pc_1^{(n)}\to\Pc_1^{(n+1)}$.
The colimit is the quotient of the tame product
$\prod_{i\geq 1}^{\text{tame}}\, \Pc_1$ by the sum 
$\bigoplus_{i\geq 1}\, \Pc_1$; so $\pi_0$ of the product is non-zero and
even has a non-trivial $M$-action.
\end{eg}

\begin{eg}[Loop and suspension]\label{eg-loop and suspension}
The loop $\Omega X$ and suspension $S^1\sm X$ of a symmetric spectrum~$X$
are defined by applying the functors $\Omega$ respectively $S^1\sm-$
levelwise, where the structure maps do not interact with the new loop 
or suspension coordinates. We claim that loop and suspension
simply shift the homotopy groups while leaving the $M$-action unchanged.

We use the isomorphism 
$\alpha:\pi_{k+n}\Omega(X_{n})\iso\pi_{1+k+n}X_{n}$
defined by sending a representing continuous map $f:S^{k+n}\to\Omega (X_{n})$
to the class of the adjoint $\hat f:S^{1+k+n}\to X_{n}$
given by $\hat f(s\sm t)=f(t)(s)$, where $s\in S^1$, $t\in S^{k+n}$.
As $n$ varies, these particular isomorphisms are compatible with
the symmetric group actions and stabilization maps, so they
form an isomorphism of $I$-functors 
$\alpha:\un{\pi}_k(\Omega X)\iso\un{\pi}_{1+k}X$, 
hence induce an isomorphism of $M$-modules 
$$ \alpha \ : \ \pi_k(\Omega X)\ \xra{\ \iso\ } \ \pi_{1+k}X $$
on colimits. 

For every symmetric spectrum $X$ 
the map $S^1\sm -:\pi_{k+n}X_n\to\pi_{1+k+n}(S^1\sm X_n)$
is $\Sigma_n$-equivariant and a natural transformations of $I$-functors
as $\mathbf n$ varies. So it induces a natural map 
$$\beta\ :\ \pi_{k}X\ \to\ \pi_{1+k}(S^1\sm X) $$ 
which is $M$-linear and an isomorphism by `stable excision`,
see~\cite[Lemma~3.1.13]{HSS}. 

Moreover, the composite
$$ \pi_* X \ \xra{\ \beta\ } \ \pi_{1+*}(S^1\sm  X) 
\ \xra{\ \alpha^{-1}} \  \pi_{*} (\Omega(S^1\sm X)) $$ 
is the map induced by the adjunction unit $X\to \Omega(S^1\sm X)$ on homotopy.
\end{eg}

\begin{eg}[Shift]\label{eg-shift}
The {\em shift} is another construction for symmetric spectra which
reindexes the homotopy groups, but unlike the suspension, this construction
changes the $M$-action in a systematic way.
The shift of a symmetric spectrum $X$ is given by
$$ (\sh X)_n \ = \ X_{1+n} $$
with action of $\Sigma_n$ via the monomorphism 
$(1\times -):\Sigma_n\to\Sigma_{1+n}$ which is explicitly
given by $(1\times \gamma)(1)=1$ and $(1\times \gamma)(i)=\gamma(i-1)+1$ 
for $2\leq i\leq 1+n$.
The structure maps of $\sh X$ are the reindexed structure maps for $X$.
An isomorphic description of the shift is as a 
symmetric function spectrum~\cite[Def.~2.2.9]{HSS}
$$ \sh X\ =\ \Hom_S(F_1S^0,X) \ ,$$
which shows that the shift functor has a left adjoint given
by $Y\mapsto F_1S^0\sm Y$.

If we view $\Sigma_n$ as the subgroup of $M$ of maps which fix all
numbers bigger than $n$, then the homomorphism 
$(1\times -):\Sigma_n\to\Sigma_{1+n}$ has a natural extension to a
monomorphism $(1\times -):M\to M$ given by $(1\times f)(1)=1$ 
and $(1\times f)(i)=f(i-1)+1$ for $i\geq 2$. 
The image of the monomorphism $1\times -$ is 
the submonoid of those $g\in M$ with $g(1)=1$. 
If $W$ is an $M$-module, we denote by $W{(1)}$ the $M$-module
with the same underlying abelian group, but with $M$-action
through the endomorphism $1\times -$. We call $W(1)$ the {\em shift} of~$W$.
Since $|1\times f|=1+|f|$, shifting an $M$-module shifts
the filtration subgroups, i.e., we have $M(1)^{(n)}=M^{(1+n)}$
for all $n\geq 0$.
Thus the $M$-module  $W{(1)}$ is tame if and only if $W$ is.

For any symmetric spectrum $X$, integer $k$ and large enough $n$
we have 
$$\pi_{(k+1)+n}(\sh X)_n \ = \ \pi_{k+(1+n)} X_{1+n} \ , $$
and the maps in the colimit system for $\pi_{k+1}(\sh X)$
are the same as the maps in the colimit system for $\pi_{k} X$.
Thus we get  $\pi_{k+1}(\sh X)=\pi_{k}X$ as abelian groups.
However, the action of a permutation on $\pi_{k+1+n}(\sh X)_n$
is shifted by the homomorphism $1\times -$, so we have 
\begin{equation}\label{homotopy of shift} 
\pi_{*+1}(\sh X)\ = \ (\pi_{*}X)(1) 
\end{equation}
as $M$-modules. 

Shifting preserves $\pi_*$-isomorphisms because of~\eqref{homotopy of shift},
but shifting does {\em not} in general
preserve stable equivalences of symmetric spectra.
An example is the fundamental stable equivalence $\lambda:F_1S^1\to S$
of~\cite[Example~3.1.10]{HSS} which is adjoint to the identity of~$S^1$.
The symmetric spectrum $\sh(F_1S^1)$ is isomorphic to the
wedge of $F_0S^1$ and $F_1S^2$, while $\sh S\iso F_0 S^1$;
the map $\sh\lambda:\sh(F_1S^1)\to\sh S$ is the projection to a wedge
summand whose complementary summand $F_1S^2$ is not stably contractible.
\end{eg}

\begin{eg}[Shift adjoint]\label{eg-drift} 
We calculate the effect on homotopy groups of the left adjoint of shifting 
by establishing a natural isomorphism of $M$-modules
\begin{equation}\label{eq-drift iso} 
\pi_{k}(F_1S^0\sm X) \ \iso \ 
\mZ[M]^+\tensor_{\mZ[M]} \pi_{k+1}X \ . 
\end{equation}
Here $\mZ[M]^+$ denotes the monoid ring of $M$ with
its usual left action, but with right action through the
monomorphism $(1\times -):M\to M$ given by $(1\times f)(1)=1$ 
and $(1\times f)(i)=f(i-1)+1$ for $i\geq 2$. 
As a right $M$-module, $\mZ[M]^+$ is free of countably infinite rank
(one possible basis is given by the transpositions
$(1,n)$ for $n\geq 1$). So the isomorphism~\eqref{eq-drift iso}
in particular implies that the underlying abelian group
of $\pi_{k}(F_1S^0\sm X)$ is a countably infinite sum of copies of the
underlying abelian group of $\pi_{k+1}X$.

The spectrum $F_1S^0\sm Y$ is trivial in level~0 and in positive levels we have
$$ (F_1S^0\sm X)_{1+n} \ = \ \Sigma_{1+n}^+ \sm_{\Sigma_{n}} X_{n}
\ . $$
Here $\Sigma_n$ acts from the right on $\Sigma_{1+n}$ via the monomorphism
$(1\times-):\Sigma_n\to\Sigma_{1+n}$.
The structure map
$(\Sigma_{1+n}^+ \sm_{\Sigma_{n}} X_{n})\sm S^1\to
\Sigma_{1+n+1}^+ \sm_{\Sigma_{n+1}} X_{n+1}$ is induced by
$(-\times 1):\Sigma_{1+n}\to\Sigma_{1+n+1}$ (the `inclusion`)
and the structure map of $X$.
Taking stable homotopy groups we get a $\Sigma_{1+n}$-equivariant
isomorphism
\begin{align*} \pi_{k+(1+n)}^{\text{s}}(F_1S^0\sm X)_{1+n} \ &= \ 
\pi_{k+1+n}^{\text{s}}\left(\Sigma_{1+n}^+ 
\sm_{\Sigma_{n}} X_{n}\right)\\
& \iso\ \mZ[\Sigma_{1+n}]\tensor_{\mZ[\Sigma_{n}]}
\pi_{(k+1)+n}^{\text{s}}X_n \ .
\end{align*} 
Taking the colimit over the stabilization maps
gives an isomorphism
$$ \pi_{k}(F_1S^0\sm X) \ \iso \ 
\mZ[\Sigma_{\infty}]^+\tensor_{\mZ[\Sigma_{\infty}]}\pi_{k+1}X $$
where $\Sigma_{\infty}$ is the subgroup of $M$ consisting of bijections
which fix almost all elements of~$\omega$.
The isomorphism~\eqref{eq-drift iso} is then obtained from the observation
that for every tame $M$-module $W$, the natural map
$\mZ[\Sigma_{\infty}]^+\tensor_{\mZ[\Sigma_{\infty}]}W
\to \mZ[M]^+\tensor_{\mZ[M]}W$ is a bijection.

We note that the functor $\mZ[M]^+\tensor_{\mZ[M]} -$ 
is left adjoint to $\Hom_M(\mZ[M]^+,-)$,
which is a fancy way of writing the algebraic shift functor
$W\mapsto W(1)$.
Under the isomorphism~\eqref{eq-drift iso} and the 
identification~\eqref{homotopy of shift}, the adjunction between
$\sh=\Hom_S(F_1S^0,-)$ and $F_1S^0\sm-$ as functors of symmetric spectra
corresponds exactly to the adjunction between
$W\mapsto W(1)$ and $\mZ[M]^+\tensor_{\mZ[M]} -$ 
as functors of tame $M$-modules.
\end{eg}

\begin{eg}[Homotopy of $R^\infty X$]\label{eg-R infty} 
In the proof of~\cite[Thm.~3.1.11]{HSS}
Hovey, Shipley and Smith define spectra $RX$ and $R^\infty X$ which
come up when studying the relationship between $\pi_*$-isomorphisms
and stable equivalences and which reappear in the characterization
of semistable symmetric spectra in~\cite[Prop.~5.6.2]{HSS}. 
We exhibit a natural isomorphism of $M$-modules
\begin{equation}\label{eq-homotopy of R^infty}
\pi_k(R^\infty X) \ \iso \ (\pi_kX)(\infty) \ . 
\end{equation}
Here for an $M$-module $V$ we denote by $V(\infty)$ the colimit
of the sequence
$$ V \ \xra{\ d\cdot\ }\ V(1) \ \xra{\ d\cdot\ } \ V(2) \ \xra{\ d\cdot\ }
\cdots$$
(note that $(1\times f)d=df$ for all $f\in M$, 
which means that $d\cdot :V\to V{(1)}$ is indeed $M$-linear, 
and so the colimit $V(\infty)$ is naturally an $M$-module.)

By definition we have $RX=\Hom_S(F_1S^1,X)=\Omega(\sh X)$.
Let $\lambda: F_1S^1\to S$ denote the morphism 
adjoint to the identity of $S^1$; this is the prototype of a stable equivalence
which is not a $\pi_*$-isomorphism, compare~\cite[Example~3.1.10]{HSS}.
Then $\lambda$ induces a morphism on function spectra
$$ \lambda^* \ : \ X=\Hom_S(S,X)\ \to \ \Hom_S(F_1S^1,X)=RX $$
and $R^\infty X$ is the colimit of the sequence 
$$ X \ \xra{\ \lambda^*}\ RX \ \xra{\ R(\lambda^*)\ }\ R^2X \ 
\xra{\ R^2(\lambda^*)\ }\ \cdots \ . $$
For calculating the homotopy of $R^\infty X$ we identify the
effect of $\lambda^*:X\to RX=\Omega(\sh X)$ in homotopy.
The level $n$ component $\lambda^*_n:X_n\to (RX)_n=\Omega( X_{1+n})$ 
is adjoint to the composite
$$ S^1\sm X_n \xra[\ \text{twist}\ ]{\ \iso\ } X_n\sm S^1 \xra{\ \sigma_n\ }
X_{n+1} \xra{(1,\dots,n,n+1)} X_{1+n}  $$
(using only the structure map $\sigma_n$ without the 
twist isomorphism and cycle permutation 
$(1,\dots,n,n+1)$ does not yield a morphism of symmetric spectra !)
So the square
$$\xymatrix@C=25mm{ \pi_{k+n}X_n \ar[d]_{\pi_{k+n}(\lambda_n^*)} 
\ar[r]^-{\iota_*} & 
\pi_{k+n+1} X_{n+1} \ar[d]^{(-1)^{n}(1,\dots,n,n+1)_*}\\
\pi_{k+n}\Omega( X_{1+n}) \ar[r]_{(-1)^k\alpha}^{\iso} &
\pi_{k+1+n} X_{1+n} }
$$
commutes, where $\iota_*$ is the stabilization map,
and the isomorphism $\alpha$ is as in Example~\ref{eg-loop and suspension}.
The signs arise as the effect 
of moving a sphere coordinate past~$k$ respectively~$n$ other coordinates.
As $n$ increases, the maps 
$(-1)^{n}(1,\dots,n,n+1)_*\circ\iota_*:\pi_{k+n}X_n\to \pi_{k+1+n} X_{1+n}$
stabilize to the left multiplication of $d\in M$ on $\pi_kX$, 
see Example~\ref{eg-d acts}. So we have shown that the square
\begin{equation}\label{eq-d is lambda}
\xymatrix@C=25mm{ \pi_{k}X \ar[d]_-{\pi_k(\lambda^*)} 
\ar[r]^-{d\cdot} & (\pi_k X)(1) \ar@{=}[d] \\
\pi_{k} RX \ar[r]_-{(-1)^k\alpha}^{\iso}  & \pi_{k+1}(\sh X) }
\end{equation}
commutes. If we iterate applications of $R$ and pass to the colimit,
we obtain the isomorphism~\eqref{eq-homotopy of R^infty}.
\end{eg}

\begin{rk}\label{rk-all operations} 
The monoid $M$ gives essentially all
natural operations on the homotopy groups of symmetric spectra.
More precisely, we now identify the ring of natural 
operations $\pi_0X\to\pi_0X$ with a completion of the monoid ring $\mZ[M]$.
Moreover, tame $M$-modules can equivalently be described
as the discrete modules over the ring of operations. 
We will not need this information later, so we will be brief.

We define the ring $\mZ[[M]]$ as the endomorphism ring of the
functor $\pi_0:\spec\to\Ab$. So an element of $\mZ[[M]]$
is a natural self-transformation of the functor $\pi_0$,
and composition of transformations gives the product. 
The following calculation of this ring depends on the fact
that the homotopy group functor $\pi_0$ is pro-represented, in the
level homotopy category of symmetric spectra, by the inverse system
of free symmetric spectra $F_nS^n$, and that we know $\pi_0(F_nS^n)$
by Example~\ref{eg-free}. 

In more detail:
for every $n\geq 0$ we let $j_n\in\pi_n (F_nS^n)_n$ be 
the wedge summand inclusion $S^n\to \Sigma_n^+\sm S^n=(F_nS^n)_n$
indexed by the unit element of~$\Sigma_n$.
Then evaluation at $j_n$ is a bijection
$$ [F_nS^n,X] \ \to \ \pi_n X_n \ , \quad [f]\mapsto f_*(j_n) $$
where the left hand side means homotopy classes of morphisms
of symmetric spectra.
We write $\lambda:F_{n+1}S^{n+1}\to F_{n}S^{n}$ for the morphism adjoint
the wedge summand inclusion
$S^{n+1}\to \Sigma_{n+1}^+\sm (S^n\sm S^1)=(F_nS^n)_{n+1}$
indexed by the unit element of~$\Sigma_{n+1}$.
Then we have
$$  \lambda_*(j_{n+1})\ =\ \iota_*(j_{n}) $$
in the group $\pi_{n+1}(F_nS^n)_{n+1}$ which implies that the squares
$$\xymatrix@C=10mm{  [F_nS^n,X] \ar[r]^-\iso\ar[d]_{[\lambda,X]} 
& \pi_n X_n\ar[d]^{\iota_*}\\
[F_{n+1}S^{n+1},X] \ar[r]_-\iso& \pi_{n+1} X_{n+1}}$$
commute. 
Passage to colimits give a natural isomorphism
$$ \colim_n\ [F_nS^n,X] \ \to \ \pi_0X \ . $$
From here the Yoneda lemma shows that 
we get an isomorphism of abelian groups
\begin{equation}\label{Z[[M]] and pi_0} 
\beta \ : \ \mZ[[M]] \ \to \ \lim_n\   \pi_0(F_nS^n) \ , 
\end{equation}
(where the limit is taken over the maps $\pi_0\lambda$)
by sending a natural transformation \mbox{$\tau:\pi_0\to\pi_0$}
to the tuple $\{\tau_{F_nS^n}[j_n]\}_{n}$.

It remains to exhibit the ring $\mZ[[M]]$ as a completion of the
monoid ring $\mZ[M]$. The natural action of
$M$ on the 0th homotopy group of a symmetric spectrum 
provides a ring homomorphism $\mZ[M]\to\mZ[[M]]$.
We define a left ideal $I_n$ of $\mZ[M]$ as the
subgroup generated by all differences  of the form $f-g$ 
for all $f,g\in M$ such that $f$ and $g$ agree on $\mathbf n$.
If $W$ is a tame $M$-module and if $x\in W^{(n)}$ has filtration~$n$,
then $I_n\cdot x=0$. So the action of the monoid ring 
$\mZ[M]$ on any tame module automatically extends to an additive map
$$ (\lim_n\ \mZ[M]/I_n) \tensor W \ \to \ W \ . $$
(Warning: $I_n$ is {\em not} a right ideal for $n\geq 1$, so
the completion does not a priori have a ring structure).
Since the homotopy groups of every symmetric spectrum form tame
$M$-modules, this gives a map of abelian groups
$$ \alpha\ :\ \lim_n\ \mZ[M]/I_n\ \to\ \mZ[[M]] $$ 
which extends the map from the monoid ring $\mZ[M]$. 

To prove that $\alpha$ is a bijection we show that the composite
$\beta\alpha:\lim_n\mZ[M]/I_n\to\lim_n\pi_0(F_nS^n)$
with the isomorphism~\eqref{Z[[M]] and pi_0} is bijective.
But this holds because the composite arise from compatible isomorphisms
$$ \mZ[M]/I_n\ \to \ \pi_0(F_nS^n)\ , \quad 
f+I_n \ \longmapsto \ f\cdot[j_n] \ , $$
which in turn uses the isomorphisms $\Pc_n\iso\pi_0(F_nS^n)$ 
from Example~\ref{eg-free}.

We end this remark by claiming without proof that the extended action
of $\mZ[[M]]$ on a tame $M$-module $W$ makes it a {\em discrete module}
in the sense that the action map
$$ \mZ[[M]] \times W \ \to \ W$$
is continuous with respect to the discrete topology on $W$
and the filtration topology on $\mZ[[M]]$.
Conversely, if $W$ is discrete module over $\mZ[[M]]$, 
then its underlying $M$-module is tame.
This establishes an isomorphism between the category of
tame $M$-modules and the category of discrete $\mZ[[M]]$-modules.
\end{rk}

\section{Semistable symmetric spectra}\label{sec-semistable}

The semistable spectra form an important class of symmetric spectra
since between these, stable equivalences coincide with
$\pi_*$-isomorphisms. The reason behind this is that for
semistable spectra, the naively defined homotopy groups
of~\eqref{colimit system pi_k} coincide with
the `true` homotopy groups, i.e., morphisms in the
stable homotopy category from the sphere spectra.

More formally, a symmetric spectrum is {\em semistable} 
in the sense of~\cite[Def.~5.6.1]{HSS}
if some (hence any) stably fibrant replacement is a $\pi_*$-isomorphism. 
Here a {\em stably fibrant replacement} is a stable equivalence
$X\to LX$ with target an $\Omega$-spectrum. 
Many symmetric spectra which arise naturally are semistable
(compare Example~\ref{eg-semistable since stabilizing}),
and Section~5.6 of~\cite{HSS} gives some criteria for checking
semistability. We now provide a criterion
for semistability in terms of the $M$-action on the homotopy groups of
a symmetric spectra.

\begin{prop}\label{prop-semistable}
A symmetric spectrum is semistable if and only if the $M$-action on
all of its homotopy groups is trivial.
\end{prop}
\begin{proof}
By~\cite[Prop.~5.6.2]{HSS}, $X$ is semistable if and only if
the map $\lambda^*:X\to RX=\Omega(\sh X)$ 
is a $\pi_*$-isomorphism (since we work with topological spaces and not
simplicial sets, we do not need any level fibrant replacement).
By~\eqref{eq-d is lambda} this happens if and only if left multiplication 
by $d$ is an isomorphism on all homotopy groups, 
which by Lemma~\ref{lemma-algebraic consequences}~(iii)
is equivalent to a trivial $M$-action.  
\end{proof}

The `trivial $M$-action` criterion is often handy for deciding
about semistability and for showing that semistability
is preserved by certain constructions. We give a few examples of this:

\begin{eg}\label{eg-semistable since stabilizing}
By Example~\ref{eg-stabilizing} any symmetric spectrum 
whose homotopy groups stabilize has trivial $M$-action
and is thus semistable (compare Proposition~5.6.4~(2) of~\cite{HSS}).
This applies in particular to suspension spectra, $\Omega$-spectra
(possibly only from some later point on), 
or symmetric spectra associated to $\Gamma$-spaces.
By Example~\ref{M acts trivially on orthogonal spectra},
the underlying symmetric spectrum of every orthogonal
spectrum has trivial $M$-action on homotopy groups, and is thus semistable.
\end{eg}

\begin{eg}
Let $X$ be a symmetric spectrum whose homotopy groups
are dimensionwise finitely generated as abelian groups. 
Tameness then forces the $M$-action to be trivial 
(Lemma~\ref{lemma-algebraic consequences}~(iv)) and so
$X$ is semistable.
This is a strengthening of Proposition~5.6.4~(1) of~\cite{HSS},
where is it is proved that spectra with {\em finite} homotopy groups
are semistable.
\end{eg}

\begin{eg}
Example~\ref{ex-products} shows that an infinite product of
symmetric spectra with trivial homotopy groups can have homotopy groups
with non-trivial $M$-action. In particular, infinite products
of semistable symmetric spectra need not be semistable.
\end{eg}

\begin{eg}
If $f:X\to Y$ is any morphism of symmetric spectra, then the
homotopy groups of the spectra $X$, $Y$ and the mapping cone 
$C(f)=[0,1]^+\sm X\cup_fY$
are related by a long exact sequence of tame $M$-modules
(see~\cite[Thm~7.4~(vi)]{MMSS}; we use that 
the $M$-action does not change under loop and suspension).
Trivial tame $M$-modules are closed under taking submodules,
quotient modules and extensions 
(Lemma~\ref{lemma-not nec tame properties}~(iv)); 
so if two out of three graded $M$-modules
$\pi_*X$, $\pi_*Y$ and $\pi_*C(f)$ have trivial
$M$-action, then so does the third. 
Thus the mapping cone of any morphism between semistable 
symmetric spectra is semistable.

If $f:X\to Y$ is an h-cofibration~\cite[Sec.~5]{MMSS} of symmetric spectra,
or simply an injective morphism when in the simplicial context of~\cite{HSS},
then the mapping cone $C(f)$ is $\pi_*$-isomorphic to the quotient $Y/X$.
Thus if two of the spectra $X$, $Y$ and $Y/X$ are semistable, 
then so is the third.
\end{eg}

\begin{eg}
Semistability is preserved under suspension, loop, wedges, shift
and sequential colimits along h-cofibrations
(or injective morphisms when in the simplicial context)
since these operations preserve the property of
$M$ acting trivially on homotopy groups.
\end{eg}

\begin{eg}\label{eg-CW smash semistable}
For a symmetric spectrum $X$ and a pointed space $K$ we let
$K\sm X$ be the symmetric spectrum obtained by smashing $K$ levelwise
with $X$. So $(K\sm X)_n=K\sm X_n$, with $\Sigma_n$-action 
by the given action on $X_n$, and with structure map
$\Id\sm\sigma_n:K\sm X_n\sm S^1\to K\sm X_{n+1}$. For example,
when $K=S^1$ is the circle, this specializes to the suspension of $X$.
We claim that if $X$ is semistable and $K$ is a CW-complex, then
the symmetric spectrum $K\sm X$ is again semistable.

We first prove the claim for finite dimensional CW-complexes by
induction over the dimension. If $K$ is 0-dimensional, then
$K\sm X$ is a wedge of copies of $X$, thus semistable.
If $K$ has positive dimension~$n$ and $K_{(n-1)}$ is its $(n-1)$-skeleton,
then $K/K_{(n-1)}$ is a wedge of $n$-spheres and so 
the quotient of $K\sm X$ by the subspectrum $K_{(n-1)}\sm X$
is a wedge of $n$-fold suspension of $X$. 
By induction the subspectrum  $K_{(n-1)}\sm X$ is semistable; 
since the inclusion is an h-cofibration and the quotient spectrum 
is also semistable, so is $K\sm X$.
For a general CW-complex $K$ the symmetric spectrum $K\sm X$ 
is the sequential colimit, over h-cofibrations,
of the smash product of $X$ with the skeleta of $K$. 
So $K\sm X$ is semistable.

The geometric realization of any simplicial set is a CW-complex,
so in the simplicial context of~\cite{HSS} we deduce that
for any pointed simplicial set $K$ and any semistable symmetric spectrum
$X$ the symmetric spectrum $K\sm X$ is again semistable.
\end{eg}

\begin{eg}\label{eg-semistable homotopy colimit}
Let $F:J\to\spec$ be a functor from a small category $J$ to
the category of symmetric spectra. If $F(j)$ is semistable for
each object $j$ of $J$, then the homotopy colimit of $F$ over $J$
is semistable. 

Indeed, the homotopy colimit is the geometric realization
of the {\em simplicial replacement} $\amalg_*F$ in the sense of Bousfield and 
Kan~\cite[Ch.~XII, 5.1]{BK}, a simplicial object of symmetric
spectra. The spectrum of $n$-simplices of $\amalg_*F$ is a wedge, 
indexed over the $n$-simplices of the nerve of $J$, of spectra
which occur as values of $F$. The geometric realization
$|\amalg_*F|$ is the sequential colimit, over h-cofibrations,
of the realizations of the skeleta $\sk_n\amalg_*F$ in 
the simplicial direction, so it suffices to show that each of these
is semistable. The skeleton inclusion realizes to an
h-cofibration $|\sk_{n-1}\amalg_*F|\to |\sk_n\amalg_*F|$
whose quotient symmetric spectrum is a wedge, 
indexed over the {\em non-degenerate} $n$-simplices of the nerve of $J$, 
of $n$-fold suspensions of spectra which occur as values of $F$. 
So the quotient spectra are semistable, 
and so by induction the symmetric spectra
$|\sk_n\amalg_*F|$ are semistable.
\end{eg}

As the final result in this section we show 
that the smash product of two semistable symmetric spectra 
is again semistable, under a mild cofibrancy condition.
We switch to symmetric spectra of simplicial set
(as opposed to topological spaces) for the rest of this section, 
because several results
which we will want to quote are in the literature in this context.
A symmetric spectrum of simplicial sets
$X$ is {\em S-cofibrant} if for all $n\geq 0$
the natural map $L_nX\to X_n$ is injective.
Here $L_nX=(X\sm \bar S)_n$ is the $n$th latching space in the
sense of~\cite[Def.~5.2.1]{HSS} (this characterization of $S$-cofibrations
does not seem to be in the literature, but the proof is
analogous to the proof of~\cite[Prop.~5.2.2]{HSS} for stable cofibrations).
A symmetric spectrum of simplicial sets
$X$ is {\em stably cofibrant} if for all $n\geq 0$
the natural map $L_nX\to X_n$ is injective and the symmetric group $\Sigma_n$ 
acts freely away from the image.
Thus every stably cofibrant symmetric spectrum is $S$-cofibrant,
but not vice versa. For example, semifree symmetric spectra
$H_nL$ associated to a pointed $\Sigma_n$-simplicial set $L$ are $S$-cofibrant,
but not stably cofibrant, unless the $\Sigma_n$-action is
free away from the basepoint.

\begin{prop}\label{smash preserves pi_*-isos} 
Smashing with an $S$-cofibrant symmetric spectrum preserves
$\pi_*$-isomorphisms.
\end{prop}
\begin{proof} We go through a sequence of steps proving that
successively larger classes of $S$-cofibrant spectra have the desired property.
By Example~\ref{eg-loop and suspension} respectively 
Example~\ref{eg-drift}, the homotopy groups of $S^1\sm X$
and of $F_1S^0\sm X$ are functors of the homotopy groups of $X$.
So smashing with  $S^1$ and $F_1S^0$ preserves $\pi_*$-isomorphisms.
Since $F_nS^m$ is isomorphic to $(F_1S^0)^{\sm m}\sm (S^1)^{\sm n}$,
smashing with the cofibrant spectrum $F_nS^m$ preserves
$\pi_*$-isomorphisms. Since homotopy groups of a wedge are
the direct sum of the homotopy groups, smashing with a wedge of
spectra of the form $F_nS^m$, for varying $n$ and $m$, preserves
$\pi_*$-isomorphisms.

Every stably cofibrant symmetric spectrum is a retract of
one built via the small object argument as  
$X=\colim X_n$, starting with $X_0=*$, and such that 
each $X_i\to X_{i+1}$ is a stable cofibration with quotient isomorphic
to a wedge of symmetric spectra of the form $F_nS^m$.
Inductively, using that cofibre sequences of spectra
give rise to long exact sequences of homotopy groups,
smashing with each $X_i$ preserves $\pi_*$-isomorphisms.
Thus smashing with the filtered colimit $X$ and thus with
an arbitrary stably cofibrant spectrum, preserves $\pi_*$-isomorphisms.

Finally, let $X$ be an $S$-cofibrant symmetric spectrum.
We choose a level equivalence $X'\to X$ with stably cofibrant source.
If $f:A\to B$ is a $\pi_*$-isomorphism, then in the commutative square
$$\xymatrix{ X'\sm A \ar[r]\ar[d] & X'\sm B \ar[d]\\
X\sm A\ar[r] & X\sm B} $$
the upper horizontal map is a $\pi_*$-isomorphism by the above.
Both vertical maps are level equivalences by~\cite[Sec.~5]{HSS}.
Thus the lower map is a $\pi_*$-isomorphism, which finishes the proof.
\end{proof}

\begin{prop}\label{prop-smash preserves semistable} 
Let $X$ and $Y$ be two semistable spectra one of which is
$S$-cofibrant. Then the smash product $X\sm Y$ is semistable.
\end{prop}
\begin{proof} Suppose $X$ is $S$-cofibrant and semistable.
We first prove the proposition when $Y$ has a special form, namely
$Y=\Omega^n L'(\Sigma^\infty K)$ for a pointed
simplicial set $K$, where $L'$ is a level fibrant replacement functor. 
Smashing with an $S$-cofibrant spectrum preserves level 
equivalences~\cite[Thm.~5.3.7]{HSS},
so $X\sm L'(\Sigma^\infty K)$  is level equivalent to $X\sm \Sigma^\infty K$,
which is isomorphic to the symmetric spectrum
$K\sm X$ and thus semistable by Example~\ref{eg-CW smash semistable}.

The counit of the adjunction
between loop and suspension is a $\pi_*$-isomorphism
$\varepsilon:S^n\sm \Omega^n L'(\Sigma^\infty K)\to L'(\Sigma^\infty K)$,
so by Proposition~\ref{smash preserves pi_*-isos}
the map 
$$ \Id\sm\varepsilon\ :\ X\sm S^n\sm \Omega^n L'(\Sigma^\infty K)
\ \to\ X\sm L'(\Sigma^\infty K)$$ 
is a $\pi_*$-isomorphism.
Since the target is semistable, so is 
$X\sm S^n\sm \Omega^n L'(\Sigma^\infty K)$. A symmetric spectrum
is semistable if and only if its suspension is, so we 
conclude that $X\sm \Omega^n L'(\Sigma^\infty K)$ is semistable.

To prove the general case we use Shipley's detection 
functor~\cite[Sec.~3]{shipley-THH}
which associates to every symmetric spectrum $Y$ the semistable
symmetric spectrum $DY$. Here $DY$ is the homotopy colimit of a functor 
$\mathcal D_Y:I\to \spec$ from the category $I$ 
to the category of symmetric spectra with values
$\mathcal D_Y(\mathbf n)=\Omega^nL'(\Sigma^\infty Y_n)$. 
By the above the symmetric spectrum 
$X\sm \mathcal D_Y(\mathbf n)=X\sm \Omega^nL'(\Sigma^\infty Y_n)$
is semistable  for each $n\geq 0$. 
Hence the homotopy colimit of the functor 
$X\sm \mathcal D_Y:I\to \spec$, which is isomorphic to $X\sm DY$,
is semistable by Example~\ref{eg-semistable homotopy colimit}.

By~\cite[Cor.~3.1.7]{shipley-THH} the semistable spectrum $Y$ is related
by a chain of $\pi_*$-isomorphisms to the symmetric spectrum $DY$.
By Proposition~\ref{smash preserves pi_*-isos},
$X\sm Y$ is thus related by a chain of $\pi_*$-isomorphisms 
to the symmetric spectrum $X\sm DY$, which we just recognized as semistable.
Hence $X\sm Y$ is semistable, which finishes the proof.
\end{proof}

\section{Shipley's spectral sequence for true homotopy groups}
\label{sec-spectral sequence}

By the `true` homotopy groups of a symmetric spectrum $X$ we mean the 
morphisms in the stable homotopy category from sphere spectra to~$X$.
In practice, one is usually interested in the true homotopy groups 
of a symmetric spectrum and not the naively defined and sometimes pathological
homotopy groups as in~\eqref{colimit system pi_k}.
However, the naive homotopy groups are more readily computable
from an explicit presentation of the symmetric spectrum, so one would
like a way to obtain the true homotopy groups from the naive ones.

The true homotopy groups can be calculated as the naive homotopy groups 
of the target of any {\em stably fibrant replacement}, i.e.,
any stable equivalence $X\to LX$ whose target is an
$\Omega$-spectrum. However, stably fibrant replacements are typically 
obtained by the small object argument, so it is often
difficult to get one's hands on their homotopy groups.
In the paper~\cite{shipley-THH} Shipley introduces a {\em detection functor}
$D$ such that for every symmetric spectrum $X$, the naive homotopy
groups of $DX$ are naturally isomorphic to the true homotopy groups of $X$.
Moreover, the spectrum $DX$ is defined as a homotopy colimit, so it
comes with a Bousfield-Kan spectral sequence for calculating 
its homotopy groups. We now discuss the spectral sequence 
associated to $DX$ from the perspective of the $M$-action
and reinterpret its $E^2$-term as Tor groups over the monoid ring of $M$. 
Then we discuss several examples
in which one can completely understand the spectral sequence.

We recall Shipley's construction from~\cite[Sec.~3]{shipley-THH},
which is inspired by B{\"o}kstedt's use of a homotopy colimit
over the category $I$ in the original definition of topological Hochschild
homology~\cite{boekstedt-THH}. 
First, Shipley
associates to a symmetric spectrum $X$ a functor 
$\mathcal D_X:I\to \spec$ from the category $I$ of finite sets and
injections to the category of symmetric spectra by the rule
$$ \mathcal D_X({\mathbf n}) \ = \ \Omega^n (F_0X_n) \ . $$
(Shipley works simplicially and also needs a level fibrant replacement
before taking loops, which is unnecessary in the topological context).
Here $F_0X_n$ is the free symmetric spectrum generated in level~0 
by the pointed space $X_n$ (so $F_0X_n$ is another name for the symmetric
suspension spectrum $\Sigma^{\infty}X_n$).

This assignment becomes a functor in $\mathbf n\in I$ as follows.
The symmetric group $\Sigma_n=I(\mathbf n,\mathbf n)$
acts on $\Omega^n (F_0X_n)=\map(S^n,F_0X_n)$
by conjugation, using the given action on $X_n$ and the
coordinate permutations on the source sphere $S^n$.
The morphism 
$\mathcal D_X(\iota):\mathcal D_X({\mathbf n})\to\mathcal D_X({\mathbf {n+m}})$
induced by the inclusion $\iota:\mathbf n\to\mathbf{n+m}$ is given in level~$k$
as follows. If $f:S^n\to X_n\sm S^k$ is an element of
$\mathcal D_X({\mathbf n})_k=\Omega^n(X_n\sm S^k)$, then 
$\mathcal D_X(\iota)(f)$ is the composite
$$ S^{n+m} \ \xra{f\sm \Id}  X_n\sm S^k \sm S^m \ \xra{\Id\sm \tau} \
 X_n\sm S^m \sm S^k \ \xra{\ \sigma^m\sm\Id\ } X_{n+m}\sm S^k \ . $$
This uniquely extends to a functor on the category~$I$.

Shipley then defines $DX$ as the homotopy
colimit of $\mathcal D_X$ over the category $I$
and proves in~\cite[Thm.~3.1.6]{shipley-THH}
that $DX$ is related by a natural chain of
$\pi_*$-isomorphisms to a stably fibrant replacement of~$X$.
So the naive homotopy groups of $DX$ are isomorphic to the
true homotopy groups of~$X$.

Since the naive homotopy groups are a homology theory on symmetric
spectra, the homotopy colimit $DX$ comes with a Bousfield-Kan
spectral sequence~\cite[Ch.~XII, \S 5]{BK} 
converging to the homotopy groups of $DX$
whose $E^2_{p,q}$-term equals $\colim_I^p \pi_q \mathcal D_X$,
the $p$th left derived functor of the colimit functor, applied to the
$I$-functor $\mathbf n\mapsto \pi_q \mathcal D_X(\mathbf n)$. 
We have
$$ \pi_q \mathcal D_X(\mathbf{n}) \ = \ 
\pi_q(\Omega^n F_0X_n) \ \iso \ \pi_{q+n}^{\text{s}}X_n \ ,  $$
and as $\mathbf n$ varies, this yields an isomorphism of $I$-functors
$\pi_q \mathcal D_X\iso\un{\pi}_{q}^{\text{s}}X$.
We show in Proposition~\ref{derived colim iso Tor} below
that there are natural isomorphisms of abelian groups
$$ \colim_I^p\ \left(\un{\pi}_{q}^{\text{s}}X \right)  \ 
\iso \ \Tor^{\mZ[M]}_p(\mZ,\pi_q X) \ . $$
This makes it clear that the derived colimits only depend on
the colimit of the functor $\un{\pi}_{q}^{\text{s}}X$,
i.e., the homotopy groups of $X$,
together with the $M$-action.

Using these isomorphisms, the Bousfield-Kan spectral sequence
takes the form of a strongly convergent half-plane spectral sequence 
\begin{equation} \label{naive2true} 
E^2_{p,q} = \Tor_p^{\mZ[M]}(\mZ,\pi_q X) \ \Longrightarrow \   
\pi_{p+q}(DX) \ . 
\end{equation}
The spectral sequence is natural in~$X$ with 
$d^r$-differential of bidegree $(-r,r-1)$.
The edge homomorphism 
$$ \mZ\tensor_M \pi_qX = E^2_{0,q} \ \to \ \pi_q (DX) $$
is induced on homotopy groups by the stably fibrant replacement $X\to LX$ 
composed with Shipley's isomorphism~\cite[Thm.~3.1.6~(1)]{shipley-THH} 
between the homotopy groups of $DX$ and those of
the stably fibrant replacement $LX$.
Here, and below, we use the notation $-\tensor_M-$ as short hand 
for the tensor product over the monoid ring $\mZ[M]$.

We will see below that the spectral sequence~\eqref{naive2true} 
collapses in many cases, for example for semistable
symmetric spectra and for free symmetric spectra
(see Example~\ref{eg-collapse free}), and it always 
collapses rationally (see Example~\ref{rational}). 
The spectral sequence typically does not collapse for
semifree symmetric spectra, see Example~\ref{naive2true for semifree}.

For the identification of derived colimits with Tor groups
we need two preparatory Lemmas. I owe the proof of Lemma~\ref{BM contractible} 
to Neil Strickland.

\begin{lemma}\label{BM contractible}
The classifying space $BM$ of the monoid $M$ is contractible.
\end{lemma}
\begin{proof} The classifying space $BM$ is the geometric
realization of the nerve of the category $\underline B M$
with one object whose monoid of endomorphisms is $M$.
Let $t\in M$ be given by $t(i)=2i$. We define an injective endomorphism
$c_t:M\to M$ as follows. For $f\in M$ and $i\in\omega$ we set
$$ c_t(f)(i) \ =\ \begin{cases}
\quad i & \text{\ if $i$ is odd, and}\\
2\cdot f(i/2) & \text{\ if $i$ is even. }
\end{cases}$$
Even though $t$ is not bijective, the endomorphism $c_t$ behaves
like conjugation by $t$ in the sense that 
the formula $c_t(f)\cdot t=t\cdot f$ holds.
Thus $t$ provides a natural transformation from the identity functor
of $\underline B M$ to $\underline B(c_t)$.
On the other hand, if $s\in M$ is given by $s(i)=2i-1$, then
$c_t(f)\cdot s=s$ for all $f\in M$, so $s$ provides a natural transformation 
from the constant functor of $\underline B M$ with values $1\in M$
to $\underline B(c_t)$. Thus via the homotopies induced by $t$ and $s$,
the identity of $BM$ is homotopic to a constant map, so $BM$
is contractible.
\end{proof}

\begin{lemma}\label{Z[M] plus and vanishing Tor} 
{\em (i)} Let $\mZ[M]^+$ denote the monoid ring of $M$ with
its usual left action, but with right action through the
monomorphism $(1\times -):M\to M$ given by $(1\times f)(1)=1$ 
and $(1\times f)(i)=f(i-1)+1$ for $i\geq 2$. 
Then for every $n\geq 0$ the map
\begin{align*}\kappa \ : \ \Pc_{1+n} &\to \ \mZ[M]^+\tensor_M \Pc_{n} 
\end{align*}
which sends the generator $(1,\dots,n+1)$ to the element
$1\tensor(1,\dots,n)$ of filtration $n+1$ in $\mZ[M]^+\tensor_M\Pc_n$ 
is an isomorphism of $M$-modules.

{\em (ii)} For every $n\geq 0$ and every abelian group $A$, the groups 
$\Tor^{\mZ[M]}_p(\mZ,\Pc_n\tensor A)$ vanish in positive dimensions.
\end{lemma}
\begin{proof}
(i) For any $n$-tuple $(x_1,\dots,x_n)$ of pairwise distinct natural numbers
we can choose $g\in M$ with $g(i)=x_i$ for $1\leq i\leq n$.
Because of $$ f\tensor(x_1,\dots, x_n) \ = \ f\tensor g(1,\dots,n) \ = \ 
f(1\times g)\cdot (1\tensor (1,\dots,n)) $$
the element $1\tensor (1,\dots,n)$ generates
$\mZ[M]^+\tensor_M \Pc_{n}$, so the map $\kappa$ is surjective.
The map $\mZ[M]^+\tensor_M \Pc_{n}\to\Pc_{1+n}$
which sends $f\tensor(x_1,\dots, x_n)$ to $(f(1),f(x_1+1),\dots,f(x_n+1))$
is right inverse to $\kappa$ since the composite sends the generator 
$(1,\dots,n+1)$ to itself. So $\kappa$ is also injective.

(ii) The groups $\Tor^{\mZ[M]}_p(\mZ,A)$ 
are isomorphic to the singular homology groups with
coefficients in $A$ of the classifying space $BM$ of the monoid $M$.
This classifying space is contractible by Lemma~\ref{BM contractible},
so the groups $\Tor^{\mZ[M]}_p(\mZ,A)$ vanish for $p\geq 1$,
which proves the case~$n=0$.

For $n\geq 1$, the $M$-modules $\Pc_{1+n}\tensor A$ and 
$\mZ[M]^+\tensor_M\Pc_{n}\tensor A$ are isomorphic by part~(i).
Since the $M$-bimodule $\mZ[M]^+$ is free as a left and right module
separately, the balancing property of Tor groups yields
\begin{align*}
\Tor^{\mZ[M]}_*(\mZ,\Pc_{1+n}&\tensor A)\ \iso \ 
\Tor^{\mZ[M]}_*(\mZ,\mZ[M]^+\tensor_M\Pc_{n}\tensor A)\\
&\iso \ \Tor^{\mZ[M]}_*(\mZ\tensor_M\mZ[M]^+,\Pc_{n}\tensor A)\
\iso \ \Tor^{\mZ[M]}_*(\mZ,\Pc_{n}\tensor A) 
\end{align*}
since $\mZ\tensor_M\mZ[M]^+$ is again the trivial right $M$-module $\mZ$.
So induction on $n$ shows that the groups
$\Tor^{\mZ[M]}_p(\mZ,\Pc_n\tensor A)$ vanish in positive dimensions.
\end{proof}

\begin{prop} \label{derived colim iso Tor}
For every $I$-functor $F$ there are natural isomorphisms of abelian groups 
$$ \colim^p_I F \ \iso \ \Tor^{\mZ[M]}_p(\mZ,\colim_\mN F) \ . $$
for all $p\geq 0$.
\end{prop}
\begin{proof} We show that the collection of functors
$$ \left\lbrace 
F\ \mapsto\  \Tor^{\mZ[M]}_p(\mZ,\colim_{\mN}F) \right\rbrace_{p\geq 0} $$ 
has all the properties that characterize the left derived functors 
of the colimit. Clearly
$$ \colim_I F \ \iso \ \mZ\tensor_M (\colim_{\mN}F) \ ,$$
so the functors agree for $p=0$. 
Taking colimit over the filtered category $\mN$ is exact,
so the Tor functors of the colimit take short exact sequences of $I$-functors
to long exact sequences of abelian groups.
The least obvious part is that the Tor groups 
$\Tor^{\mZ[M]}_p(\mZ,\colim_{\mN}F)$ 
vanish for all projective functors $F$ and all $p\geq 1$.
The $I$-functors $\mZ[I(\mathbf n,-)]$ 
arising as the linearizations of the representable functors
form a set of projective generators for the category of $I$-functors,
so it suffices to show the vanishing of higher Tor groups for these.
But the colimit over $\mN$ of the $I$-functor $\mZ[I(\mathbf n,-)]$
is precisely the $M$-module $\Pc_n$;
so Lemma~\ref{Z[M] plus and vanishing Tor}~(ii) provides the vanishing 
result and finishes the proof.
\end{proof}

\begin{eg}[Semistable and free symmetric spectra]\label{eg-collapse free}
When $X$ is a semistable symmetric spectrum
or a free symmetric spectrum, then the higher Tor groups
for the homotopy of $X$ vanish by Lemma~\ref{Z[M] plus and vanishing Tor}~(ii).
Thus in the spectral sequence~\eqref{naive2true}
we have $E^2_{p,q}=0$ for $p\neq 0$, and so the edge homomorphism
$$ \mZ\tensor_M(\pi_*X) \ \to \pi_*(DX) $$
is an isomorphism.
\end{eg}

\begin{eg}[Eilenberg-Mac Lane spectra] 
In Example~\ref{ex-EM ring spectra with M}
we associate an Eilenberg-Mac Lane spectrum $HW$ to every
tame $M$-module $W$. The homotopy groups of $HW$ are concentrated
in dimension~0, where we get the module $W$ back.
So the spectral sequence~\eqref{naive2true} for $HW$ collapses 
onto the axis $q=0$ to isomorphisms
$$  \pi_p(D(HW)) \ \iso \ \Tor^{\mZ[M]}_p(\mZ,W) \ . $$
In particular, the true homotopy groups of $HW$ need not be concentrated
in dimension~0.

Here is an example which shows that for non-trivial $W$ the  Eilenberg-Mac Lane
spectrum $HW$ can be stably contractible: we let $W$ be the kernel of
a surjection $\Pc_n\to\mZ$. 
Lemma~\ref{Z[M] plus and vanishing Tor} and the long exact sequence 
of Tor groups
show that the groups  $\Tor^{\mZ[M]}_p(\mZ,W)$ vanish for all $p\geq 0$.
Thus the homotopy groups of $D(HW)$ are trivial, 
i.e., $HW$ is stably contractible.
\end{eg}

\begin{eg}[Rational collapse]\label{rational} 
We claim that for every tame $M$-module $W$ and all $p\geq 1$,
we have $\Tor^{\mZ[M]}_p(\mQ,W)=0$.
So the spectral sequence~\eqref{naive2true} always collapses rationally 
and the edge homomorphism is a rational isomorphism
$$ \mQ\tensor_M(\pi_*X) \ \to \ \mQ\tensor \pi_* (DX) \ . $$
The rational vanishing of higher Tor groups is special for
{\em tame} $M$-modules.

To prove the claim we consider a monomorphism $i:V\to W$ of tame $M$-modules
and show that the kernel of the map 
$\mZ\tensor_Mi:\mZ\tensor_MV\to \mZ\tensor_MW$
is a torsion group. 
The inclusions $W^{(n)}\to W$ induce an isomorphism
$$ \colim_n\ \mZ\tensor_{\Sigma_n} W^{(n)} \ \xra{\ \iso\ } \ 
\mZ\tensor_MW\ . $$
For every $n\geq 0$, the kernel of 
$\mZ\tensor_{\Sigma_n} i^{(n)}:\mZ\tensor_{\Sigma_n} V^{(n)}
\to \mZ\tensor_{\Sigma_n} W^{(n)}$
is annihilated by the order of the group $\Sigma_n$. 
Since the kernel of $\mZ\tensor_Mi$ is the colimit of the kernels
of the maps $\mZ\tensor_{\Sigma_n} i^{(n)}$, it is torsion.
Thus the functor $\mQ\tensor_M-$ 
is exact on short exact sequences of tame $M$-modules and
the higher Tor groups vanish as claimed.
\end{eg}

The Tor groups discussed in the next lemma arise in the spectral 
sequence~\eqref{naive2true} for semifree symmetric spectra.

\begin{lemma}\label{lemma-derived of induced}
For every $\Sigma_n$-module $B$ we have a natural isomorphism
$$  \Tor^{\mZ[M]}_*(\mZ,\Pc_n\tensor_{\Sigma_n}B) \ \iso \ 
H_*(\Sigma_n;B) \ . $$
\end{lemma}
\begin{proof} Since $\Pc_n$ is free as a right $\Sigma_n$-module, the functor
 $\Pc_n\tensor_{\Sigma_n}-$ is exact.
The functor takes the free $\Sigma_n$-module of rank~1 to $\Pc_n$,
so by Lemma~\ref{Z[M] plus and vanishing Tor}~(ii)
it takes projective $\Sigma_n$-modules 
to tame $M$-modules which are acyclic for the functor $\mZ\tensor_M-$. 

Thus if $P_\bullet\to B$ is a projective resolution of 
$B$ by $\Sigma_n$-modules,
then $\Pc_n\tensor_{\Sigma_n}P_\bullet$ is a resolution
of $\Pc_n\tensor_{\Sigma_n}B$ which can be used to calculate the
desired Tor groups. Thus we have isomorphisms
$$ \Tor^{\mZ[M]}_*(\mZ,\Pc_n\tensor_{\Sigma_n}B) \ =\ 
H_*(\mZ\tensor_M\Pc_n\tensor_{\Sigma_n}P_\bullet)\ \iso \
H_*(\mZ\tensor_{\Sigma_n}P_\bullet)\ = \ H_*(\Sigma_n;B)\ . $$
\end{proof}

\begin{eg}[Semifree symmetric spectra] \label{naive2true for semifree}
For semifree symmetric spectra (see Example~\ref{ex-semifree}) 
the spectral sequence~\eqref{naive2true} typically does not degenerate. 
As an example we consider the semifree symmetric spectrum
$H_2S^2$, where $S^2$ is a $\Sigma_2$-space by coordinate permutations.

We first identify the stable equivalence type of $H_2S^2$.
The spectrum $H_2S^2$ is isomorphic to the quotient spectrum
of $\Sigma_2$ permuting the smash factors of $(F_1S^1)^{\sm 2}$.
Since the $\Sigma_2$-action on $(F_1S^1)^{\sm 2}$ is free,
the map 
$$(F_1S^1)^{\sm 2}\sm_{\Sigma_2}E\Sigma_2^+\ 
\to\ (F_1S^1)^{\sm 2}/\Sigma_2=H_2S^2$$
which collapses $E\Sigma_2$ to a point is a level equivalence. 
On the other hand, the  stable equivalence 
$\lambda^{\sm 2}:(F_1S^1)^{\sm 2}\to S$ 
is $\Sigma_2$-equivariant, so it induces a stable equivalence
$$(F_1S^1)^{\sm 2}\sm_{\Sigma_2}E\Sigma_2^+\ 
\to\ S\sm_{\Sigma_2}E\Sigma_2^+ = \Sigma^\infty B\Sigma_2^+ $$
on homotopy orbit spectra. Altogether we conclude that $H_2S^2$ is stably
equivalent to $\Sigma^\infty B\Sigma_2^+$. 

The spectral sequence~\eqref{naive2true} for $H_2S^2$ has as $E^2$-term
the Tor groups of $\pi_*(H_2S^2)$. 
According to~\eqref{homotopy of semifree spectrum} these homotopy groups
are isomorphic to~$\Pc_2\tensor_{\Sigma_2}(\pi_{*+2}^{\text{s}}S^2)(\sgn)$.
The sign representation cancels the sign action induced
by the coordinate flip of $S^2$,
so we have an isomorphism of $M$-modules 
$\pi_q(H_2S^2)\iso\Pc_2\tensor_{\Sigma_2}\pi_{q}^{\text{s}}S^0$,
this time with trivial action on the stable homotopy groups of spheres.
Using Lemma~\ref{lemma-derived of induced}, the 
spectral sequence~\eqref{naive2true} for $H_2S^2$ takes the form
$$ E^2_{p,q} \ \iso \ H_p(\Sigma_2;\pi_q^{\text{s}}S^0) \ \Longrightarrow \
\pi_{p+q}^{\text{s}}  (B\Sigma_2^+) \ . $$
This spectral sequence has non-trivial differentials and it seems likely
that it coincides with the Atiyah-Hirzebruch spectral sequence for the
stable homotopy of the space $B\Sigma_2^+$.
\end{eg}

\end{document}